\documentclass[10pt]{article} 

\usepackage[utf8]{inputenc} 

\usepackage{geometry} 
\usepackage{graphicx} 
\usepackage{xcolor}
\usepackage{parskip}
\usepackage{bm}
\usepackage{amsmath}
\usepackage{amssymb}
\usepackage{mathrsfs}
\usepackage{subfig}
\usepackage{footnpag}
\usepackage[round]{natbib}
\usepackage{indentfirst}

\usepackage{bm}
\usepackage{subfig}
\usepackage{footnpag}

\usepackage{arydshln}
\usepackage{hyperref}

\usepackage{listings}
\usepackage{doi} 

\usepackage{algorithm,algpseudocode}

\newenvironment{breakablealgorithm}
  {
   \begin{center}
     \refstepcounter{algorithm}
     \hrule height.8pt depth0pt \kern2pt
     \renewcommand{\caption}[2][\relax]{
       {\raggedright\textbf{\fname@algorithm~\thealgorithm} ##2\par}%
       \ifx\relax##1\relax 
         \addcontentsline{loa}{algorithm}{\protect\numberline{\thealgorithm}##2}%
       \else 
         \addcontentsline{loa}{algorithm}{\protect\numberline{\thealgorithm}##1}%
       \fi
       \kern2pt\hrule\kern2pt
     }
  }{
     \kern2pt\hrule\relax
   \end{center}
  }

\title{Rapid nonlinear convex guidance using a monomial method}
\author{%
Ethan R. Burnett\thanks{Marie Skłodowska-Curie Postdoctoral Fellow, Department of Aerospace Science and Technology, Politecnico di Milano, 20156 Milan, Italy \texttt{ethanryan.burnett@polimi.it}} \ and Francesco Topputo\thanks{Professor, Department of Aerospace Science and Technology, Politecnico di Milano, 20156 Milan, Italy}}

\date{} 
        
\begin{document}
\maketitle
\begin{abstract} 
This paper addresses the challenge of accommodating nonlinear dynamics and constraints in rapid trajectory optimization, envisioned for use in the context of onboard guidance. We present a novel framework that uniquely employs overparameterized monomial coordinates and pre-computed fundamental solution expansions to facilitate rapid optimization while minimizing real-time computational requirements. The fundamental solution expansions are pre-computed using differential algebra. Unlike traditional approaches that repeatedly evaluate the nonlinear dynamics and constraints as part of complex shooting or collocation-based schemes, this method replaces the nonlinearity inherent to dynamics and constraint functions entirely with a computationally simpler manifold constraint. With this approach, trajectory optimization is posed efficiently as a path planning problem on the manifold. This problem is entirely convex except for the manifold constraint, readily lending itself to solution via sequential convex programming. We demonstrate the effectiveness of our approach in computing fast and accurate delta-V optimal solutions for long-range spacecraft rendezvous, including problems with nonlinear state constraints.
\end{abstract}

\section{Introduction}
In the context of modern guidance, navigation, and control (GNC), spacecraft autonomy presents itself as a particularly technically and politically challenging problem. Strict computational and hardware limitations are imposed by the need for radiation-hardened processors and on-board power and mass constraints. Experimentation with autonomous agents is furthermore extremely regulated in comparison to terrestrial ventures in autonomy (such as self-driving cars) due to the high risk of testing autonomous capabilities with no flight heritage on extremely expensive space vehicles. In practice, this has always motivated onboard guidance implementations that are computationally lean, deterministic, and easy to test and validate. Nonetheless recent trends in dropping cost-to-orbit and the proliferation of CubeSats anticipate the deployment of numerous lower-cost deep space spacecraft that will require an increased degree of autonomy to avoid overwhelming human-in-the-loop on-ground tracking and planning capabilities \citep{DiDomenico2023}. Thus it seems likely that there will be increasing demand for on-board guidance, including for use in critical planning and decision-making operations, in the years to come. Such ``computational guidance and control" schemes \citep{CGC}, making use of greater embedded computation capability, will go far beyond the traditional algebraic operations needed for evaluating closed-form guidance, while retaining the flexibility and trust necessary for onboard use.

The main characteristics of suitable on-board spacecraft guidance, are, according to \cite{Starek2016}, that 1) it should be computationally reasonable, 2) it should compute an optimal solution wherever possible, and 3) it should enable verifiability. Convex optimization-based guidance presents an appealing candidate for meeting these challenges of onboard guidance because by nature it is fast, extremely stable, computationally well-posed, and thus easy to profile computationally \citep{Boyd_Vandenberghe_2004}. Within the convex optimization framework, it is also possible to develop \textit{stochastic guidance} schemes that robustly accommodate expected dispersions in control performance as well as state and parameter estimation errors. Chance-constrained methods are becoming popular for this because they provide performance guarantees with a user-defined confidence level -- see e.g. \cite{Oguri_JGCD_convex, berning2024BO}. In addition to stochasticity, nonlinearity poses a perennial and ubiquitous challenge. Even for otherwise convex problems, the general non-convexity of nonlinear dynamics and constraints often forces an iterative approach via sequential convexification \citep{WangGrant_sCVX}, whereby a non-convex problem is solved locally as a convex sub-problem subject to trust region constraints to enforce a stable iterative march towards the optimal solution. While extremely powerful, sequential convexification implementations tend to require fairly problem-specific details of choice of transcription and trust-region updates to ensure stability, feasibility, and computational efficiency. 

Among the aforementioned challenges of on-board guidance, for this work we are focused on methods for fast, robust, and easily characterizable onboard guidance for nonlinear systems. Guidance of nonlinear systems invariably comes with some non-trivial computational cost. However, there is operational flexibility in the timing of these computations, namely whether they happen in \textit{real-time} (i.e. when the guidance algorithm is actively updating the planned trajectory) or at \textit{non-critical times} (i.e. when the spacecraft guidance computer has no time-sensitive tasks or is comparatively idle). Thus not all of the computational cost needs to be paid in real-time, and much can be shifted to non-critical times. Some computations can be performed on-ground days in advance, or even loaded into the guidance computer before launch.   
This is done by enabling pre-computation of useful and reusable information to the maximum extent possible. 
Along these motivational lines, we introduce a new framework for using nonlinear expansions of problem dynamics that is demonstrably computationally efficient. It requires no real-time integration, either explicitly (as in predictor-corrector methods or multiple shooting) or implicitly (as with collocation schemes). It also allows for solutions whose accuracy can be easily characterized. 

Many past works are relevant to the ideas and techniques expounded in this work. First, to facilitate our computations, we make use of a computerized monomial algebra, whereby monomials of an arbitrary order and number of variables are represented by arrays of their coefficients, whose spatial arrangement matches the ordering of the monomial terms. Similar schemes have been implemented before, and \cite{AlgebraicManip} gives a broad overview of commonly used methods. \cite{JorbaNormalForms} develops a computer algebra system in which monomials are computed, stored, and manipulated in a manner very similar to ours. Their applications are otherwise quite different, leveraging manipulation of monomials in a sequence of canonical transformations for the purpose of constructing normal forms and obtaining approximate integrals of Hamiltonian systems.

The other half of our methodology involves the computation of nonlinear fundamental solution expansions about a reference -- each associated with a particular unique monomial in the initial conditions. We've explored their computation by many means. Firstly, via ``State Transition Tensors" (STTs), which are an extension of the ubiquitous state transition matrix involved in almost any method involving both linearization and discretization. \cite{ParkScheeresSTTs} provide a classic introduction, and \cite{BooneMcMahon} provide a more recent predictor-corrector guidance implementation for two-burn maneuvers. We have also computed the nonlinear fundamental solutions leveraging Differential Algebra (DA) methods, which facilitates automated computation of derivatives via a structure allowing direct treatment of many topics related to the differentiation and integration of functions \citep{Berz_ParticleBeamMaps}. The resulting representation of the expansion is called a Taylor map \citep{Berz_ParticleBeamMaps}. This methodology has seen noteworthy use in spaceflight GNC in uncertainty propagation \citep{Valli_DA_2013} and optimal control \citep{DiLizia_DA_2014,Greco_DA}. DA is our preferred method due to its versatility and speed avantages in comparison to STTs. Lastly, nonlinear fundamental solution expansions can be computed analytically for some specialized problems. This is typically aided nowadays by the use of computational software such as Mathematica, and usually involves the use of perturbation methods \citep{nayfeh2000perturbation,Hinch_Book}. In this work, our example application is one for which many such expansions have been analytically derived: the spacecraft relative motion and rendezvous problem. Works such as \cite{Butcher2016_AAS, ButcherBurnett2017} and \cite{Willis2ndOrder, willis2019secondorder} are thus especially relevant to us.

This work makes use of convex optimization for trajectory optimization. In addition to the thorough foundational work of \cite{Boyd_Vandenberghe_2004}, see also work outlining the successive convexification algorithm SCvx \citep{sCVX1,sCVX2}, which allows optimization of non-convex nonlinear systems via an iterative approach, with some guarantees of convergence not often found elsewhere in literature. This successive convexification is also occasionally referred to as a ``sequential convex programming" (SCP) method \citep{sCVX2}, but we note here its distinction from the more general sequential implementations in literature \citep{HofmannTopputo_JGCD2021, Morelli_2022}. Recent improvements and augmentations of this and SCP methods are numerous (e.g. \cite{Oguri_sCVX} applies an augmented Lagrangian formulation to enable a feasibility guarantee of SCvx). Python is our chosen language for convex optimization prototyping, offering mature open-source tools such as CVXPY \citep{diamond2016cvxpy,agrawal2018rewriting} and ECOS \citep{ECOS_2013}, as well as the recent introduction of CVXPygen \citep{cvxpygen2022} for especially rapid problem-solving via generation of a speedy custom solver implementation in C called directly from the Python implementation as a CVXPY solver method. For a clear and simple introduction to convex optimization methods for spacecraft trajectory optimization, see \cite{WangGrant_sCVX}, and for a discussion of the typically necessary collocation schemes, we recommend \cite{Kelly_colloc}. \cite{Betts_Optimization} provides a classic overview of optimization methods for trajectory design, and the recent survey papers \cite{Malyuta_OptimizationSurvey, Wang_ConvexOptimizationSurvey} showcase the strength and popularity of modern convex optimization tools and techniques for aerospace vehicle guidance and control. 

Our work follows some notable other works making use of convex optimization for spacecraft trajectory optimization. Much of this work falls under the ERC-funded EXTREMA project \citep{DiDomenico2023} for self-driving interplanetary CuebSats. \cite{HofmannTopputo_JGCD2021} present a computationally simple and robust convex optimization-based algorithm for low-thrust interplanetary trajectories. \cite{Morelli_2022} applies convex optimization to a similar problem while considering a homotopic energy-to-fuel optimal approach along with ``second-order" trust region methods. \cite{HofmannJSR_CVX} performs a wide study of different discretization and trust region methods for convex low-thrust trajectory optimization. Finally, \cite{Saglioano_PoweredDescent}, \cite{Szmuk_PoweredDescent}, and \cite{SaglianoAndLu_PoweredDescent} provide applications of convex optimization to another continuous-thrust problem: powered descent trajectory optimization. This problem in particular has enjoyed great success in convex optimization. 
In this early work we lay the groundwork of our methodology for both continuous and impulsive control, but we explore only the latter in-depth. Our example application is convex optimization of the long-range spacecraft rendezvous problem under impulsive thrust, and the methodology proposed is agnostic to the dynamical environment or the type of reference orbit. This is an application that historically has involved ground-based computations and often, manual operation of the ``chaser" spacecraft executing the rendezvous. Automation is not only highly desirable but necessary to accommodate the breadth of new scenarios and architectures anticipated by emerging and maturing commercial space operations \citep{WoffindenAutonomy}. For applications to the close-range problem with convex programming, see e.g. \cite{berning2024BO} (applying drift safety guarantees) or \cite{BurnettJGCD2022} (for methods applicable to any general periodic orbits).

\subsection{Contributions}
Despite great progress in the past decade, existing methods for rapid spacecraft trajectory optimization have their deficits. \cite{berning2024BO} and \cite{BurnettJGCD2022}, along with many others, are limited to close-proximity applications due to the linearity assumption in the dynamics. The state transition tensor-based guidance of \cite{BooneMcMahon}, while valid for longer range, only considers two-burn impulsive maneuvers, whereas often there are delta-V benefits by allowing the possibility for additional burns \citep{PrussingNImpulse}. In general the works using SCP rely extensively on shooting and/or collocation methods which must repeatedly evaluate nonlinear dynamics (and possibly constraint) functions in real-time. The accuracy of collocation methods is variable, with details of their implementation somewhat heuristic, prompting complex studies of various strategies for a problem of interest \citep{HofmannJSR_CVX}. The novel contribution of this work is a framework for posing nonlinear trajectory optimization problems in a convex optimization approach that does not require real-time integration or any other evaluation of nonlinear or linearized equations of motion, minimizing real-time computational effort. The only required operations are simple computation and manipulation of monomials and manageable linear algebra. Furthermore there is no reliance on collocation, interpolation, or other heuristic approximations of the continuous dynamics at all. In essence the method approximates the nonlinear dynamics in a linear form by moving the nonlinearity from dynamics and state constraints to manifold constraints, which are easier to work with. This is done via a new ``flattened" representation of the information contained in the state transition tensor/Taylor map, and a new monomial-based nonlinear state representation of ``osculating initial conditions". The resulting transcription scheme is easy to set up, and furthermore facilitates simple predictions of the accuracy of computed solutions based on the expansion order. The methodology enables rendezvous and station-keeping guidance solutions with practically any number of impulsive maneuvers to be computed. In this foundational work we consider only the case that the reference is a natural (i.e. control-free) trajectory, hence the limited application to rendezvous and station-keeping. 

\subsection{Structure}
This paper is organized as follows. In Section 2, we highlight the fundamental ideas behind our methodology. We discuss the ideas of osculating initial conditions and overparameterized monomial representations, which are both central and necessary to this work. 
In Section 3 we move on to the applications of our framework to spacecraft trajectory optimization, culminating in an example simple sequential convexification implementation whose real-time operations are computationally extremely lean and easy to implement.  
In Section 4 we provide a thorough example application to the nonlinear spacecraft rendezvous problem. Section 4.1 provides a simple two-stage (linear prediction, nonlinear correction) convex guidance strategy leveraging our developments. Finally, Section 4.2 showcases illustrates the sequential convexification scheme from Section 3, and in Sections 4.3 and 4.4 the principles explained in this paper are used to develop more complex examples with superior working coordinates and nonlinear path constraints. 
Section 5 has concluding remarks.

\section{Theoretical Developments}
\subsection{Osculating Initial Conditions}
We start the discussion of our methodology by considering the following controlled nonlinear system: 
\begin{equation}
\label{theory1}
\dot{\bm{x}} = \bm{f}(\bm{x}, \bm{u},t)
\end{equation}
where $\bm{u} \in \mathbb{R}^{m}$ is a control signal, $\bm{x} \in \mathbb{R}^{N}$ is a state deviation, and $t$ is time. Note that in our notation, vector quantities and vector-valued functions are bolded, and two-dimensional matrices are unbolded capital letters (Latin or Greek). 
This form also admits the simpler control-affine case where control manifests in the right-hand side of Eq.~\eqref{theory1} as a term like $+B(\bm{x},t)\bm{u}$, and also the sub-cases where $B(t)$ is not a function of the state, or is constant altogether. 
The methodology in this work is developed for systems where the state, subject to short limited intervals of intense control effort, is well-approximated across these intervals by discrete jumps. The application is developed for the dynamics of a space vehicle but should work for any system satisfying this property. 
In general, our methodology is nonlinear but still \textit{local} -- we consider the controlled motion of a state deviation $\bm{x}(t)$ in the vicinity of a natural 
reference trajectory $\bm{X}_{r}(t)$, or similarly, an equilibrium point $\bm{X}_{r}^{*}$. To avoid burdensome notation, we never show the $\delta$ on such state deviations from a reference, and hereafter we use $\bm{X}_{r}(t)$ when needed to refer to any reference, with $X_{r,i}$ for the $i$\textsuperscript{th} component of the reference state vector, and lowercase $\bm{x}(t)$ for any departure, with $x_{i}$ for the $i$\textsuperscript{th} component of the state deviation. Thus always ``$\delta\bm{x}$" $ \equiv \bm{x}$, whereby the origin $\bm{x} = \bm{0}$ denotes a point of zero departure from some reference/equilibrium point.

\begin{figure}[h!]
\centering
\includegraphics[]{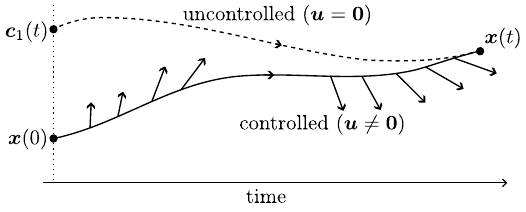}
\caption{Osculating Initial Conditions}
\label{fig:OscICs}
\end{figure}

For any valid initial condition of a deviation (henceforth simply a ``state") $\bm{x}(0)$ and a time $t>0$, we can get the state $\bm{x}(t)$ by application of the \textit{flow of the natural dynamics} (i.e., control-free, $\bm{u} = \bm{0}$) as below:
\begin{equation}
\label{theory2pre}
\bm{x}(t) = \bm{\varphi}(\bm{x}(0),t,0)
\end{equation}
where e.g. $\bm{\varphi}(\bm{x},t_{2},t_{1})$ propagates a state $\bm{x} \in \mathbb{R}^{N}$, without control, from $t_{1}$ to $t_{2}$. Similarly the inverse mapping obtains the initial conditions implied by the state at time $t$ \textit{from the flow of the natural dynamics}:
\begin{equation}
\label{theory2}
\bm{x}(0) = \bm{\varphi}^{-1}(\bm{x}(t), t, 0)
\end{equation}
Based on Eq.~\eqref{theory2}, we could also take some controlled state $\bm{x}(t)$ at some time $t>0$ and back-propagate using the flow of the natural dynamics to a corresponding ``osculating" initial condition at time $0$. We denote this osculating initial condition as $\bm{c}_{1}(t)$. This is depicted in Fig.~\ref{fig:OscICs}. As the trajectory evolves, the osculating initial condition will assume different values over time. Consider that $\bm{c}_{1}(0) = \bm{x}(0)$, but at time $0 + \delta t$, the osculating initial condition will have moved a bit if any control was applied in the interval $[0, \delta t]$. By a general time $t$, it will have moved smoothly under the action of smooth control to the depicted location, or impulsively, via jump(s), under the action of impulsive control. The mapping to and from $\bm{x}(t)$ is always well-defined for our problems of interest, so the $\bm{c}_{1}$ serves as an equally valid coordinate description of the problem. In this sense, for a controlled state $\bm{x}(t)$, there is always an ``osculating" initial state $\bm{c}_{1}(t)$: 
\begin{equation}
\label{theory2b}
\bm{c}_{1}(t) = \bm{\varphi}^{-1}(\bm{x}_{c}(t), t, 0)
\end{equation}
where the subscript ``c" emphasizes that the state $\bm{x}(t)$ may be subject to some kind of control. Conveniently, this osculating state is stationary whenever control is not being applied.

For a linear dynamical system, the transformation $\bm{\varphi}$ is linear and easy to compute, and so past works have explored solving linearized guidance problems leveraging parameterization in terms of the osculating initial conditions (or more generally, integration constants) of the system \citep{Guffanti_ICs,BurnettJGCD2022}. This transcribes trajectory optimization problems as path planning problems in a transformed domain where the state $\bm{c}_{1}(t)$ moves only via control, with otherwise null dynamics. For nonlinear systems, by contrast, we are accustomed to computing the mapping $\bm{\varphi}$ or its inverse via relatively costly numerical integration techniques. Because of this, before this work a general reformulation of nonlinear guidance or control problem in terms of initial states would not be considered much of an improvement on the problem to be solved (and indeed, would seem unnecessarily complex). In this work we develop a practical method valid for nonlinear systems in the vicinity of some reference.

\subsection{Overparameterized Representation}
Instead of a minimal parameterization (linear in the initial conditions), a complex but critical development is to instead consider the parameterization of our system by $\bm{c}_{j}$ -- a collection of \textit{all unique multivariate monomials} up to $j$\textsuperscript{th} order, generated from the $N$ (osculating) initial states $x_{1}(0), \ldots, x_{N}(0)$. Thus the aforementioned $\bm{c}_{1}$ is just $\bm{c}_{j}$ for $j = 1$. For $j > 1$, $\bm{c}_{j} \subseteq \mathbb{R}^{K_{j}}$, with $K_{j}$ given below: 
\begin{equation}
\label{WhatIsKj1}
K_{j} = \sum_{q = 1}^{j} \begin{pmatrix} N + q - 1 \\ q \end{pmatrix}
\end{equation}
One approach for building an ordered $\bm{c}_{j}$ is as follows. For all unique multivariate monomials of a given order $r$, multiply by $x_{1}(0)$, and append to the list. Then multiply the order-$r$ multivariate monomials by $x_{2}(0)$, and append in order only the non-repeated values, repeating this procedure through all state variables to $x_{N}(0)$. This process initializes with $\bm{c}_{1} = \bm{x}(0)$ and $r = 1$, and finishes after $r = j-1$ to produce $\bm{c}_{j}$. For example, for $N = 3$, $j = 2$, $\bm{c}_{2}$ is computed as below:
\begin{equation}
\label{c_build_demo}
\bm{c}_{2} = \begin{pmatrix}x_{1}(0), \ x_{2}(0), \ x_{3}(0), \ x_{1}^{2}(0), \ x_{1}x_{2}(0), \ x_{1}(0)x_{3}(0), \ x_{2}^{2}(0), \ x_{2}(0)x_{3}(0), \ x_{3}^{2}(0)\end{pmatrix}^{\top}
\end{equation}
We provide functions to perform this procedure in the Appendix A. See e.g. \cite{AlgebraicManip,JorbaNormalForms} for more details about the computerized algebraic manipulation of sets of monomials. 

With the construction of $\bm{c}_{j}$ established, we can in general nonlinearly expand the solution of a dynamical system about a fixed point (or an arbitrary reference, w.l.o.g.) in terms of unique monomials and their associated fundamental solutions as below:
\begin{equation}
\label{NLExp1}
\begin{split}
\bm{x}(t) \approx & \ x_{1}(0)\bm{\psi}_{x_{1}}(t) + x_{2}(0)\bm{\psi}_{x_{2}}(t) + \ldots + x_{N}(0)\bm{\psi}_{x_{N}}(t) \\ & \ + x_{1}^{2}(0)\bm{\psi}_{x_{1}^{2}}(t) + x_{1}(0)x_{2}(0)\bm{\psi}_{x_{1}x_{2}}(t) + \ldots +  x_{N}^{2}(0)\bm{\psi}_{x_{N}^{2}}(t) \\ & \ + \ldots \\ & \ + x_{1}^{j}(0)\bm{\psi}_{x_{1}^{j}}(t) + x_{1}^{j-1}(0)x_{2}(0)\bm{\psi}_{x_{1}^{j-1}x_{2}}(t) + \ldots + x_{N}^{j}(0)\bm{\psi}_{x_{N}^{j}}(t) \\ = & \ \Psi_{j}(t)\bm{c}_{j}
\end{split}
\end{equation}
where the ellipses include all unique monomial terms not shown explicitly. The $\bm{\psi}$ functions are the partial derivatives of the trajectory with respect to monomial functions of the (osculating) initial conditions, and may be defined as below:
\begin{equation}
\label{WhatIsPsi}
\bm{\psi}_{x_{\alpha_{1}}\ldots x_{\alpha_k}}(t) = \frac{\mathcal{R}(\alpha_{1}\ldots\alpha_{k})}{k!} \frac{\partial^{k}\bm{X}_{r}(t)}{\partial X_{r,\alpha_{1}}(0)\ldots\partial X_{r,\alpha_{k}}(0)}; \ \ t \geq 0
\end{equation}
where $\mathcal{R}(\alpha_{1}\ldots\alpha_{k})$ denotes the number of unique permutations for the list $X_{r,\alpha_{1}}\ldots X_{r,\alpha_{k}}$ and $k$ is the order of the differential. For example, the list $X_{r,1}X_{r,2}X_{r,1}$ has 3 unique permutations, so $\mathcal{R}(1,2,1) = 3$. Note $k=1$ recovers the columns of the classical state transition matrix. To be clear, a more standard multi-index notation (see e.g. \cite{NeidingerMultiIndex}) gives the following equivalences to the above definitions: 
\begin{equation}
\label{MultiIndex0}
\bm{x}(t) = \sum_{|\beta| = 1}^{j}\left(\bm{\psi}_{\beta}(t)\prod_{l=1}^{N}x_{l}^{\beta_{l}}(0)\right)
\end{equation}
\begin{subequations}
\label{MultiIndex1}
\begin{align}
\bm{\psi}_{(1, 1,\ldots, 0)}(t) = & \ \bm{\psi}_{x_{1}x_{2}}(t) \\
\bm{\psi}_{(0, 1,\ldots, 2)}(t) = & \ \bm{\psi}_{x_{2}x_{N}^{2}}(t)
\end{align}
\end{subequations}
\begin{equation}
\label{MultiIndex2}
\bm{\psi}_{\beta}(t) = \frac{1}{\beta_{1}!\beta_{2}!\ldots\beta_{N}!}\left(\left(\frac{\partial}{\partial X_{r,1}(0)}\right)^{\beta_{1}}\left(\frac{\partial}{\partial X_{r,2}(0)}\right)^{\beta_{2}}\ldots\left(\frac{\partial}{\partial X_{r,N}(0)}\right)^{\beta_{N}}\bm{X}_{r}(t)\right)
\end{equation}
where $|\beta| = \beta_{1} + \beta_{2} + \ldots + \beta_{N}$, and $\beta$ contains $N$ integers (by contrast, $\alpha$ contains $k$ integers). We often opt for the explicit symbolic form of Eqs.~\eqref{NLExp1} and \eqref{WhatIsPsi}. 

In our approach, the expansion order $j$ can be chosen based on the level of nonlinearity of the problem to be solved. Note that for locally convergent expansions, the error of the representation of a particular trajectory decreases monotonically as order $j$ is increased. See for example \cite{ButcherBurnett2017} for further discussion of error vs. expansion order. The $\bm{\psi}$ functions can be constructed to a given desired order by various means, including using the relevant terms from the associated state transition tensors (STTs) for the problem \citep{ParkScheeresSTTs,BooneMcMahon} or the Taylor map (TM) from differential algebra (DA) \citep{Berz_ParticleBeamMaps, Valli_DA_2013}, computed numerically, or perturbation solutions \citep{nayfeh2000perturbation,Hinch_Book} which are computed analytically. These topics are discussed in greater detail in the Appendices B, C, and D. An important thing to note is that, at a given order, $\Psi_{j}$ is a minimal representation (i.e. none of the information in $\Psi_{j}$ is redundant) written in a linear form. By contrast, the STTs are a redundant representation (by the equivalence of mixed partials), written in general in summation form. This is a matter of combinations vs. permutations of the differentials. Eq.~\eqref{NLExp1} reduces to the familiar expression $\bm{x}(t) = \Phi(t,0)\bm{x}(0)$, where $\Phi(t,0)$ is the state transition matrix (STM) when $j = 1$, retaining only the $N$ linear fundamental solutions $\bm{\psi}_{x_{1}}(t)$ through $\bm{\psi}_{x_{N}}(t)$ in the first row of Eq.~\eqref{NLExp1}. For the general case $j > 1$, the matrix $\Psi_{j}(t)$ is $N\times K_{j}$. Thus $\Psi_{j}$ is a wide (no. columns $\geq$ no. rows) matrix that not only trivially reduces to the STM for $j=1$, it can be thought of as a flattened (and minimal) STT that still retains a simple 2D matrix form. Like the STM, this matrix is a function of time. We sometimes write this as $\Psi_{j}(t)$, implying a continuous function, but in practice it is available only at discrete times, $\Psi_{j}(t_{i})$, and furthermore we reserve the use of more explicit notation $\Psi_{j}(t_{i},t_{0}) \equiv \Psi_{j}(t_{i})$ for epoch time $t_{0}$ when needed.

\begin{figure}[h!]
\centering
\includegraphics[]{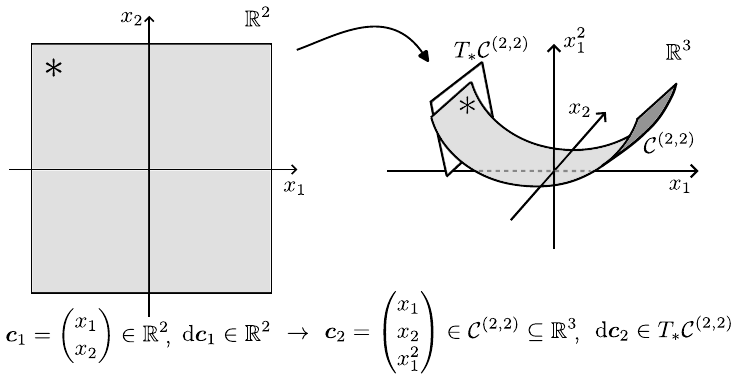}
\caption{Monomial Coordinates on a Manifold}
\label{fig:myfig1}
\end{figure}

We say that $\bm{c}_{j}$ is ``overparameterized": it has $K_{j}$ components to describe a system with $N$ states, with $K_{j} > N$ for $j>1$. Thus there are $K_{j}-N$ constraints relating its nonlinear to its $N$ linear elements. For example, applying a change in the linear component of $\bm{c}_{j}$ associated with $x_{1}(0)$, the higher-order components of $\bm{c}_{j}$ such as $x_{1}^{2}(0)$, $x_{1}(0)x_{2}(0)$, $x_{1}(0)x_{3}(0)$ etc. are forced to take on a certain value, because they are functionally dependent on the linear components. This can be shown by computing variations of Eq.~\eqref{c_build_demo} based on various values of $\bm{x}(0) \in \mathbb{R}^{N}$, and noting the consequence on the higher-order components' values. Thus $\bm{c}_{j}$ is constrained to lie on an $N$-dimensional surface that we denote as $\mathcal{C}^{(N,j)}$, which is embedded in a $K_{j}$-dimensional space. Figure~\ref{fig:myfig1} depicts for a system $\bm{x} = (x_{1}, x_{2})^{\top}$ (i.e. $N = 2$) the surface obtained with a set of monomials up to order 2 - omitting $x_{1}^{2}$ and $x_{1}x_{2}$ from the second-order monomials (reducing $K_{j}$ from 5 to 3) to enable a simple 3D drawing.

It can be shown by use of the monomial equations that the surface $\mathcal{C}^{(N,j)}$ is an $N$-dimensional sub-manifold of the Euclidean space $\mathbb{R}^{K_{j}}$ (e.g. by applying the Definition 3.10 of \cite{boumal2023intromanifolds}). Furthermore any admissible $\bm{c}_{j}$ is one-to-one with a unique $\bm{c}_{1}$ and vice-versa, thus $\mathcal{C}^{(N,j)}$ is homeomorphic to $\mathbb{R}^{N}$. This can be visualized geometrically in the simple 2D case of Figure~\ref{fig:myfig1}: $\mathbb{R}^{2}$ can be ``curled" upwards to produce $\mathcal{C}^{(2,2)}$. Given a point $\bm{c}_{j} \in \mathcal{C}^{(N,j)}$, any physically admissible infinitesimal variations $\text{d}\bm{c}_{j}$ must lie in the tangent space $T_{\bm{c}_{j}}\mathcal{C}^{(N,j)}$. Critically, this tangent space is easy to compute. Because the sequence of monomials is uniquely defined by the powers of the constituent state components, the Jacobian $\frac{\partial\bm{c}_{j}}{\partial\bm{c}_{1}}$ admits a simple analytic form, and infinitesimal deviations obey the following:
\begin{equation}
\label{WhatIsTangSpace}
\text{d}\bm{c}_{j} \in T_{\bm{c}_{j}}\mathcal{C}^{(N,j)}; \ \ \text{d}\bm{c}_{j} = \frac{\partial\bm{c}_{j}}{\partial\bm{c}_{1}}\text{d}\bm{c}_{1}
\end{equation}
The Jacobian at a point $\bm{c}_{j} \in \mathcal{C}^{(N,j)}$ can be shown to be a linear function of $\bm{c}_{j}$. The only nonlinear function necessary is the mapping from $\bm{c}_{1}$ to $\bm{c}_{j}$, and we denote the following convenient notation for the mapping between linear ($j=1$) and nonlinear ($j > 1$) monomial sequences:
\begin{subequations}
\label{Ej}
\begin{align}
\bm{c}_{j} = & \ \bm{E}_{j}(\bm{c}_{1}) \\
\bm{c}_{1} = & \ \bm{E}_{j}^{-1}(\bm{c}_{j}) = \begin{bmatrix} I_{N\times N} & 0_{N\times (K_{j}-N)} \end{bmatrix}\bm{c}_{j}
\end{align}
\end{subequations}
where $\bm{E}_{j}$ expands $\bm{c}_{1}$ to its nonlinear order $j$ representation and conversely $\bm{E}_{j}^{-1}$ extracts just the linear part $\bm{c}_{1}$ from $\bm{c}_{j}$. An efficient strategy for computing the function $\bm{E}_{j}$ is provided in Appendix A. The structure $\mathcal{C}^{(N,j)}$ has many favorable properties. It is 1) smooth, 2) analytically parameterized by the monomial equations, 3) admits efficient computations in the tangent space, 4) maps one-to-one with $\mathbb{R}^{N}$. These properties, along with others to be introduced, inspire a nonlinear optimization strategy exploiting this structure. 

\subsection{Kinematics in terms of the monomials}
To motivate more complex general arguments, we start with a simplified and specific example. Let $\bm{x}(t) = \left(\bm{r}^{\top}(t), \bm{v}^{\top}(t)\right)^{\top}$ denote specifically the state in Cartesian coordinates with position $\bm{r}(t)$ and velocity $\bm{v}(t)$. Thus the state at time $t_{i}$ obeys the following mapping from the monomial expansion of the initial conditions:
\begin{equation}
\label{TransX1}
\bm{x}(t_{i}) = \Psi_{j}(t_{i})\bm{c}_{j}
\end{equation}
Eq.~\eqref{TransX1} is only approximately true for $t_{i} > 0$, but it becomes more accurate as the order $j$ is increased (consider that the STM will not exactly describe the evolution of the state deviation from $\bm{x}(0)$ to $\bm{x}(t_{i})$, but a second-order expansion will do better within some region of convergence, and third order better still). We defer for now discussion of the error involved in this expression, and limit our arguments to some domain of validity still to be defined. 

\begin{figure}[h!]
\centering
\includegraphics[width=5.5in]{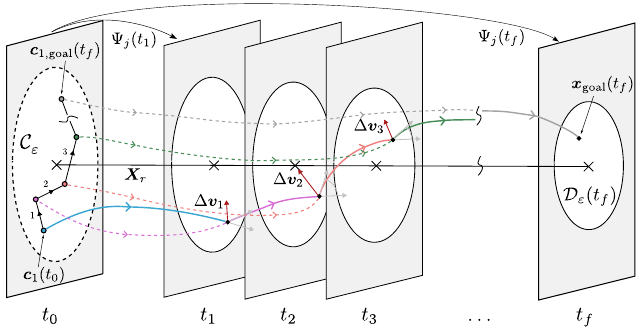}
\caption{Impulsive maneuver trajectory transcription via $\Psi_{j}$}
\label{fig:TrajTrans1}
\end{figure}

Let us furthermore limit ourselves for now to the case of control of $\bm{x}(t)$ executed via impulsive maneuvers. (This case implies, operationally, that the reference $\bm{X}_{r}(t)$ is a natural trajectory.) For this control problem, we allow for the possibility of a maneuver at any time $t_{i}$ in a discretized sequence $[t_{1}, t_{2}, t_{3}, \ldots, t_{K}]$. In lieu of considering each $\Delta\bm{v}(t_{i})$ and the product of their effects across all discretization points, we can instead parameterize the problem in terms of the true initial condition $\bm{c}_{j}(t_{0}) = \bm{E}_{j}(\bm{x}(0))$, the osculating initial condition $\bm{c}_{j,\text{goal}}(t_{f})$ which generates $\bm{x}_{\text{goal}}(t_{f})$ under application of the map $\Psi_{j}(t_{f}, t_{0})$, and a sequence of intermediate states $\bm{c}_{j}(t_{i}) = \bm{E}_{j}(\bm{c}_{1}(t_{i}))$, each of which represents an osculating initial condition that, under application of the appropriate map $\Psi_{j}(t_{i},t_{0})$, intersects with a particular maneuver node at its associated maneuver time $t_{i}$. This is illustrated by Fig.~\ref{fig:TrajTrans1}. The horizontal black line depicts the reference trajectory $\bm{X}_{r}(t)$. We depict also five trajectory nodes -- the initial time $t_{0}$, three maneuver times $t_{1}, \ t_{2}, \ t_{3}$, and the final time $t_{f}$. The true executed trajectory is depicted as a sequence of solid colored lines, and the osculating parts as dashed lines. We show explicitly the action of three maneuvers at times $t_{1}$, $t_{2}$, and $t_{3}$, and subsequent maneuvers at discretization times not depicted complete the transition of the system to the goal state. The problem is entirely parameterized in the left-most image: the domain of the osculating initial conditions. Thus each natural arc of $\bm{x}(t)$ corresponds to a single point on the left-most plane, and maneuvers induce transitions in these stationary states. Optimization of the trajectory involves a deterministic process of finding the optimal $\bm{c}_{j}(t_{i})$ (or equivalently, $\bm{c}_{1}(t_{i})$) for all $t_{i}$. In this manner the impulsive maneuver trajectory optimization problem is transformed into a path-planning problem. The white regions represent static ($\mathcal{C}_{\varepsilon}$) and evolving ($\mathcal{D}_{\varepsilon}(t_{i})$) domains of validity, still to be discussed, and the ``$\times$" mark the local origins, intercepting with the reference trajectory, i.e. where $\bm{x}(t_{i}) = \bm{0}$. 

During an impulsive maneuver applied at some $t_{i}$, the position will not change, which provides constraints on admissible jumps in $\bm{c}_{j}$. Meanwhile, the velocity change is a linear function of the jump in $\bm{c}_{j}$, under the monomial-based approximation. These describe kinematic conditions relating the instantaneous state and the overparameterized monomials of the osculating initial conditions, written explicitly below:
\begin{subequations}
\label{CoVI1}
\begin{align}
& \Delta\bm{r}(t_{i}) :=  \Psi_{r,j}(t_{i})\left( \bm{c}_{j}(t_{i}) - \bm{c}_{j}(t_{i-1}) \right) = \bm{0} \\
& \Delta\bm{v}(t_{i}) := \Psi_{v,j}(t_{i})\left( \bm{c}_{j}(t_{i}) - \bm{c}_{j}(t_{i-1}) \right)
\end{align}
\end{subequations}
where $\Psi_{r,j}$ denotes the top $N/2$ (position-associated) rows of $\Psi_{j}$, and $\Psi_{v,j}$ the bottom $N/2$ (velocity-associated) rows. This enables the complete parameterization of the trajectory in terms of decision variables $\bm{c}_{j}(t_{i}) \ \forall \ t_{i} \in [t_{1}, t_{2}, t_{3}, \ldots, t_{f}]$, which can be linearly mapped to the instantaneous state $\bm{x}(t_{i})$. We will show that this parameterization is extremely powerful. Given pre-computation of the map $\Psi_{j}$ from $t_{0}$ to all times of interest, nonlinear guidance problems can be repeatedly solved with very low recurrent computational cost, via reuse of the map $\Psi_{j}$ and manageable linear algebra. There are strong opportunities for this methodology in situations where a reasonable reference trajectory is known.

\subsection{Domain of validity of the representation}
\subsubsection{Bounding the error of the representation}
Because we use a truncated nonlinear map $\Psi_{j}$ to describe the departure of the state from the reference, in reality our methods are limited to a finite region in the vicinity of the reference trajectory. We present some arguments to enable the quantification of this region, but note we do not give a formal proof. This procedure would only need to be applied once, after the computation of the expansion. First, we assume that the dynamics in the vicinity of the reference are smooth and continuous and admit a locally convergent Taylor series expansion. We define a scalar measure of the truncation error as
\begin{equation}
\label{TruncErr1}
e_{j}(\bm{c}_{1},t_{i}) = \|\bm{x}_{\text{approx}}(t_{i}) - \bm{x}_{\text{true}}(t_{i})\| = \|\Psi_{j}(t_{i},0)\bm{E}_{j}(\bm{c}_{1}(t_{i})) - \bm{\varphi}(\bm{c}_{1}(t_{i}),t_{i},0)\|
\end{equation}
The value of $e_{j}(\bm{c}_{1},t)$ tends to increase with the norm of $\bm{c}_{1}$ (fixing $t$) and with the scale of $t$ (fixing $\bm{c}_{1}$). 

In lieu of ensuring that each temporal discretization point of a given trajectory is contained within an instantaneous domain of validity $\mathcal{D}(t_{i})$ (a local neighborhood of the reference $\bm{X}_{r}(t_{i})$ for which the expansion about the reference is sufficiently valid), we can instead ensure that the osculating initial conditions of the entire problem lie within a \textit{static} domain of validity $\mathcal{C}_{\varepsilon}$. We do this by setting an error threshold $\varepsilon$ as the max desired truncation error at time $t_{f}$, as defined by Eq.~\eqref{TruncErr1}, and we consider error sampling of a discrete set of $D$ control points $\bm{c}_{1}$ on the surface of an $(N-1)$-dimensional ball $\mathcal{B}_{\mathscr{R}}=\mathscr{R}S^{N-1}$: 
\begin{equation}
\label{whoIsCr1}
\mathcal{B}_{\mathscr{R},D} = \left\{ \bm{c}_{1}[l] \in \mathcal{B}_{\mathscr{R}}, 0 < l \leq d \right\}
\end{equation}
where $\bm{c}_{1}[l]$ denotes the $l$\textsuperscript{th} state vector in the sample set, and in other words we require $\|\bm{c}_{1}[l]\| = \mathscr{R}$ for all $0 < l \leq d$. Noting $e_{j}(\bm{0},t) :=  0$, and defining an error tolerance $\varepsilon$ for which we require $e_{j}(\bm{c}_{1},t_{f}) \leq \varepsilon$, the region of validity of a given expansion $\Psi_{j}$ can be defined as the radius satisfying the following for the given error tolerance:
\begin{equation}
\label{whoIsCr2}
\mathscr{R}_{\text{crit}}(\varepsilon) = \left\{ \mathscr{R} \in \mathbb{R}^{\geq 0} \ | \ \underset{\bm{c}_{1}\in\mathcal{B}_{\mathscr{R}}}{\text{max}}e_{j}(\bm{c}_{1},t_{f}) = \varepsilon \right\}
\end{equation}
In other words, $\mathscr{R}_{\text{crit}}(\varepsilon)$ is the minimum radius of $\mathcal{B}_{\mathscr{R}}$ such that the representation error equals $\varepsilon$ (as a maximum) somewhere on $\mathcal{B}_{\mathscr{R}}$. The maximum error on $\mathcal{B}_{\mathscr{R}}$ increases monotonically with increasing $\mathscr{R}$, so $\mathscr{R}_{\text{crit}}(\varepsilon)$ can be estimated with a bisection search approximating maximum error on $\mathcal{B}_{\mathscr{R},D}$ at various values of $\mathscr{R}$. The $\mathcal{B}_{\mathscr{R}}$ is assumed to be adequately densely sampled if increasing the value of $D$ by some factor does not appreciably change the value of the estimated maximum error.
With an estimate of $\mathscr{R}_{\text{crit}}(\varepsilon)$, the static domain of validity is then defined as: 
\begin{equation}
\label{whoIsCr3}
\mathcal{C}_{\varepsilon} := \{\bm{c}_{1} : \|\bm{c}_{1}\| \leq \mathscr{R}_{\text{crit}}(\varepsilon) \}
\end{equation}
This argument also exploits the fact that $e_{j}(\bm{c}_{1},t)$ tends to grow with increasing $t$. By bounding the error at a final time of interest $t_{f}$, the error at prior times $t_{i} < t_{f}$ should also be similarly bounded. 

\subsubsection{Certified guidance results}
The rigorous application of this method is thus summarized as follows. First, the expansion $\Psi_{j}(t)$ is computed, for a particular reference $\bm{X}_{r}$ and a particular discretization of time. Then, a numerical procedure as described previously is applied to compute $\mathscr{R}_{\text{crit}}(\varepsilon)$ and $\mathcal{C}_{\varepsilon}$, and the use of $\Psi_{j}$ is then ``certified" at error level $\varepsilon$ for any optimal control problem whose solution lies within $\mathcal{C}_{\varepsilon}$ (with a statistical level of confidence informed by the number of samples $D$). Furthermore, any guidance solution whose boundary values satisfy $\|\bm{c}_{1}(t_{0})\| < \mathscr{R}_{\text{crit}}(\varepsilon)$ and $\|\bm{c}_{1}(t_{f})\| < \mathscr{R}_{\text{crit}}(\varepsilon)$ can be constrained to lie within the region of validity via the following requirement on all trajectory nodes:
\begin{equation}
\label{CertifiedGuidance}
\|\bm{c}_{1}(t_{i})\| < \mathscr{R}_{\text{crit}}(\varepsilon) \ \forall \ i \in [1,K]
\end{equation}
This requirement, which has the property of ensuring that the guidance solution is trustworthy, can easily be augmented to the other constraints in a given optimization problem.

\subsubsection{Further considerations for estimating error}
The prior arguments make implicit assumptions about the smoothness of parameteric variations of $e_{j}(\bm{c}_{1},t_{f})$ with $\bm{c}_{1}$ -- it is possible to imagine pathological error functions whose behavior on $\mathscr{B}_{\mathscr{R}}$ confounds the sampling-based certification due to very small localized error maxima not being sampled. It is likely that more can be done to compute the region of validity efficiently. 
Bounding the error based on the order of expansion should be possible for our problems of interest -- see the truncation error arguments in DA literature, e.g. in \cite{WittigDomainSplit}. For a convergent expansion, at order $j$ the error of Eq.~\eqref{TruncErr1} will be dominated by the order $j+1$ term, so the error can be upper-bounded by $Q\|\bm{c}_{1}\|^{j+1}$ for some integer $Q$. Note lastly that the transitions in $\bm{c}_{1}$ produced by the optimization procedure (i.e. the maneuvers) will also introduce errors of order $\mathcal{O}(\|\bm{c}_{j}\|^{j+1})$ because the optimal control problem is truncated and solved at order $j$. These could have an impact on the region of validity computation. In practice, however, for our problems of interest, the order $j$ can easily be chosen sufficiently high to make the effect of these errors sub-dominant compared to other expected sources of error, particularly model and execution errors. 

In reality, the information provided by $\Psi_{j}$ for accurately re-parameterizing the whole problem at the initial time could also somehow be captured by smaller maps (i.e. lower $j$) from each $t_{i}$ to $t_{i+1}$. This is because the approximation afforded by $\Psi_{j}(t,0)$ worsens as time $t\geq 0$ increases. However, re-parameterizing the problem at the initial time accomplishes two goals: 1) it facilitates a \textit{linear} transcription of the trajectory optimization problem (extremely useful for developing our convex optimization-based schemes in Section III) and 2) it greatly simplifies the effort needed to certify (or ensure) the accuracy of the guidance solution, because we need only check (or enforce) that the optimal path lies within the \textit{static} domain of validity, e.g. $\bm{c}_{1}(t_{i}) \in \mathcal{C}_{\varepsilon} \ \forall \ t_{i}$. 

\subsection{Generalizations}
\subsubsection{Continuous theory: variation of parameters}
Here, we extend the prior discretized arguments to the case of continuous control operating in continuous time. The controlled dynamics are given by the equations below. Instead of discrete jumps, we seek the \textit{continuously} varying $\bm{c}_{j}(t)$ such that Eq.~\eqref{TransX1} continually describes the controlled state. Again we note these arguments only hold in circumstances where the truncation error is negligible. Factoring and differentiating Eq.~\eqref{TransX1}:
\begin{subequations}
\label{CoVU1}
\begin{align}
\dot{\bm{r}} = & \ \dot{\Psi}_{r,j}\bm{c}_{j} + \Psi_{r,j}\dot{\bm{c}}_{j} \\
\dot{\bm{v}} = & \ \dot{\Psi}_{v,j}\bm{c}_{j} + \Psi_{v,j}\dot{\bm{c}}_{j} 
\end{align}
\end{subequations}
Because $\Psi_{j}(t)$ dictates how $\bm{r}$ and $\bm{v}$ evolve for natural arcs, in principle higher-order state derivatives (e.g. $\bm{a} = \frac{\text{d}}{\text{d}t}(\bm{v})$) can be obtained from time derivatives of $\Psi_{j}$. The following kinematic identity must be satisfied, noting that by definition of our choice of coordinates $\bm{x}$, the nonlinear fundamental solutions satisfy $\Psi_{v,j}(t) = \dot{\Psi}_{r,j}(t)$:
\begin{equation}
\label{CoVU2}
\begin{split}
\dot{\bm{r}} = \bm{v} = & \ \Psi_{v,j}\bm{c}_{j} \\ = & \ \dot{\Psi}_{r,j}\bm{c}_{j}
\end{split}
\end{equation}
From this we obtain the below constraint on admissible directions of $\dot{\bm{c}}_{j} \in T_{\bm{c}_{j}}\mathcal{C}^{(N,j)}$:
\begin{equation}
\label{CoVU3}
\Psi_{r,j}\dot{\bm{c}}_{j} = \bm{0}
\end{equation}
Then, denoting the natural acceleration $\bm{a}$ and the control acceleration $\bm{a}_{c} = B_{L}(\bm{x},t)\bm{u}$ where $B_{L}$ is the lower $N/2$ rows of control matrix $B$ (and furthermore $B_{L} = I_{3\times 3}$ in the general 3D Cartesian case), we obtain also the following relationship from Eq.~\eqref{CoVU1} after noting $\bm{a} = \dot{\Psi}_{v,j}\bm{c}_{j}$, which itself is obtained analogously to Eq.~\eqref{CoVU2}:
\begin{equation}
\label{CoVU4}
\bm{a}_{c} = \Psi_{v,j}\dot{\bm{c}}_{j}
\end{equation}
Examining Eqs.~\eqref{CoVU3} and \eqref{CoVU4}, we recognize these expressions as a continuous control analog of the earlier kinematic constraints in Eq.~\eqref{CoVI1}. A review of variation of parameters (e.g., \cite{schaub}) might assist the interested reader. In this role the osculating initial conditions are, colloquially, analogous to the osculating orbital elements -- they are stationary in the absence of disturbance/control. Note also that in general these continuous equations can also be obtained by evaluating their discrete-time counterparts by evaluating the limit that $(t_{i}-t_{i-1})\rightarrow 0$.
\subsubsection{Functions and coordinate transformations}
Consider some general (and possibly nonlinear) scalar function of the state at a particular time, $g(\bm{x},t)$. Let $g$ be a smooth function of $\bm{x}(t_{i})$ that admits a Taylor series, thus it can be written as a multivariate polynomial in the components of $\bm{x}(t_{i})$. Because each component $x_{q}(t_{i})$ for $q = 1:N$ is itself a multivariate polynomial in the $x_{q}(0) \ \forall q$, and because multivariate polynomials are closed under multiplication and addition, the function $g(\bm{x},t)$ can also be approximated (to order $j$) as a multivariate polynomial in the $x_{q}(0)$. Thus $g(\bm{x},t)$ can itself be approximated as a \textit{linear function} of $\bm{c}_{j}$:
\begin{equation}
\label{AnyScalarFunc}
g(\bm{x}(t_{i}),t_{k}) \approx \bm{\gamma}_{j}^{\top}(t_{k})\bm{c}_{j}(t_{i})
\end{equation}
for vector $\bm{\gamma}_{j}$ of length $K_{j}$. Similarly, vector functions of the state, $\bm{g}(\bm{x},t) \in \mathbb{R}^{p}$, as collections of $p$ scalar functions, can be related linearly to $\bm{c}_{j}$:
\begin{equation}
\label{AnyScalarFunc}
\bm{g}(\bm{x}(t_{i}),t_{k}) \approx \Gamma_{j}(t_{k})\bm{c}_{j}(t_{i})
\end{equation}
for matrix $\Gamma_{j}$ of size $p \times K_{j}$. In this regard, general functions of the discretized state can be approximated by simple linear functions of the overparameterized monomial state.

There are two immediate consequence of the preceding arguments. First, any admissible nonlinear constraint functions on $\bm{x}$ that can be rendered linear in $\bm{c}_{j}$ fit within the convex optimization framework for our parameterization $\bm{c}_{j}$. Note however that the form of the original nonlinearity will still be important to practical performance of the SCP scheme. Second, the $\Psi_{j}$ and $\bm{c}_{j}$ can be developed for any desired coordinates $\bm{\eta}$, then related back to our kinematically convenient Cartesian coordinates $\bm{x}$ via use of $\bm{x} = \bm{g}(\bm{\eta},t)$. In other words, the kinematic constraints of Eq.~\eqref{CoVI1} generalize to:
\begin{subequations}
\label{CoVI1new}
\begin{align}
& \Gamma_{r,j}(t_{i})\left( \bm{c}_{j}^{(\eta)}(t_{i}) - \bm{c}_{j}^{(\eta)}(t_{i-1}) \right) = \bm{0} \\
& \Delta\bm{v}(t_{i}) = \Gamma_{v,j}(t_{i})\left( \bm{c}_{j}^{(\eta)}(t_{i}) - \bm{c}_{j}^{(\eta)}(t_{i-1}) \right)
\end{align}
\end{subequations}
where the superscript ``$(\eta)$" emphasizes that the fundamental solutions and monomial states are developed specifically in working coordinates $\bm{\eta}$. We find that the computation of $\Gamma_{j}$ is facilitated greatly by use of open-source differential algebra tools such as Pyaudi \citep{audi_izzo_zenodo}.

\section{Application to spacecraft trajectory optimization}
\begin{figure}[h!]
\centering
\includegraphics[]{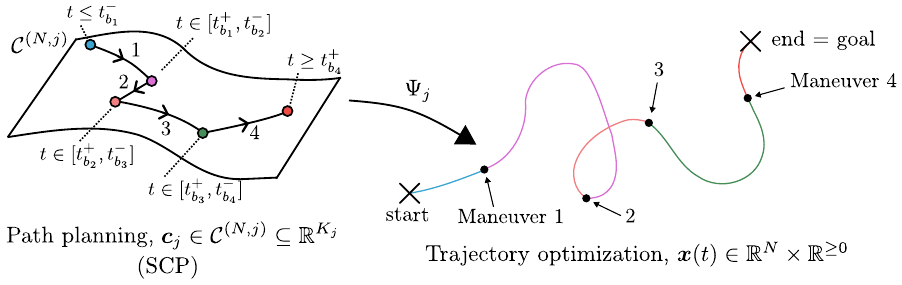}
\caption{Trajectory Optimization as Path Planning on $\mathcal{C}^{(N,j)}$, simplified}
\label{fig:PathPlans1}
\end{figure}

Consider again the optimization of nonlinear spacecraft trajectories composed of 1) impulsive maneuvers, and 2) coast arcs between maneuvers. The monomial coordinates are naturally well-suited for posing this trajectory optimization problem efficiently as a path-planning problem. We start with a simple unconstrained delta-V optimal nonlinear control example, although state constraints are explored in Section IV. 
The transcription between monomial coordinates $\bm{c}_{j} \in \mathcal{C}^{(N,j)}$ and the trajectory $\bm{x}(t)$ is shown in an example 4-burn solution in Fig.~\ref{fig:PathPlans1}, which is a simplified complement to Fig.~\ref{fig:TrajTrans1} emphasizing the non-Euclidean nature of $\mathcal{C}^{(N,j)}$. This figure depicts that the monomial state is stationary in the absence of maneuvers, and therefore any change in the monomial state indicates the application of a maneuver. In the figure the maneuver times are represented as $t_{b_{i}}$, with a ``$-$" or ``$+$" superscript indicating the instant just before or after the burn. Furthermore, the monomial state is constrained to lie on $\mathcal{C}^{(N,j)}$, but is otherwise unconstrained, with all discretization nodes free to assume their respective locations such that the overall minimizing path is obtained. 
We explore two optimization schemes: minimizing $J = \sum_{i}\|\Delta\bm{v}(t_{i})\|^{2}$ (i.e. ``energy optimal") and also $J = \sum_{i}\|\Delta\bm{v}(t_{i})\|$ (``fuel optimal"). We differentiate between them by their exponent ``$\mathcal{P}$", 2 or 1. The former is parsed explicitly in this work, whereas the objective function for the latter is implemented via CVXPY parsing.
\subsection{Problem formulation}
The unconstrained fixed-time optimization problem with impulsive maneuvers is given below:
\begin{equation}
\label{TrajOpt1}
\begin{split}
& \underset{\Delta\bm{v}(t_{i}) \ \forall \ i \in [1,K]}{\text{min}} \ \ \ \sum_{i=1}^{K}\left.\|\Delta\bm{v}(t_{i})\|\right.^{\mathcal{P}} \\ \\
& \text{subject to} \ \ \ \begin{matrix} & \ \bm{x}(t_{0}) = \bm{x}_{0} \\& \  \bm{x}(t_{K}^{+}) = \bm{x}_{f}  \\ & \ \bm{x}(t_{i+1}^{-}) = \bm{\varphi}\left(\bm{x}(t_{i}^{+}), t_{i+1}, t_{i}\right) \\ & \ \bm{r}(t_{i}^{+}) = \bm{r}(t_{i}^{-}) \\
& \ \bm{v}(t_{i}^{+}) = \bm{v}(t_{i}^{-}) + \Delta\bm{v}(t_{i})\end{matrix} 
\end{split}
\end{equation}
where we explore $\mathcal{P} = 1,2$, and both yield a convex problem, and furthermore $\bm{x}_{f} := \bm{x}_{\text{goal}}(t_{f})$ is the goal state at fixed final time $t_{f} = t_{K}$, and $\bm{x}_{0}$ is the initial condition. The trajectory is split into up to $K+1$ different points, corresponding to an initial condition plus $K$ times where maneuvers are allowed. For $K\gg N$, the optimal solution will generally have $\Delta\bm{v}(t_{i}) = \bm{0}$ for many discretized times $t_{i}$. The times are not optimization variables; time is discretized a priori. The expression $\bm{\varphi}(\bm{x}(t_{1}), t_{2}, t_{1})$ denotes the flow of the state $\bm{x}(t)$ from $t_{1}$ to $t_{2}$, and $t_{i}^{-}$ and $t_{i}^{+}$ denote the time $t_{i}$ at the instants before and after a maneuver $\Delta\bm{v}(t_{i})$. The first two conditions in Eq.~\eqref{TrajOpt1} are simple boundary conditions based on prescribed initial and final states of the trajectory. The third condition is a requirement that the state between maneuvers obey the flow of the natural dynamics. The fourth is a continuity condition requiring that the position is unchanged during an instantaneous maneuver, and the fifth records the velocity change induced by a delta-V.

We seek to rewrite the optimization problem of Eq.~\eqref{TrajOpt1} in terms of the monomial coordinates. Recall the mapping from monomial coordinates to the instantaneous state at some time $t$:
\begin{equation}
\label{TrajOpt2}
\bm{x}(t) = \Psi_{j}(t,t_{0})\bm{c}_{j}
\end{equation}
Given an instantaneous maneuver $\Delta\bm{v}(t_{i})$, the corresponding $\Delta \bm{c}_{j}(t_{i}) = \bm{c}_{j}(t_{i}) - \bm{c}_{j}(t_{i-1})$ must satisfy the following for a $\Psi_{j}$ partitioned row-wise into $\Psi_{r,j}$ (top $N/2$ rows - corresponding to position states) and $\Psi_{v,j}$ (bottom $N/2$ rows - velocity states):
\begin{subequations}
\label{TrajOpt3}
\begin{align}
\Delta\bm{c}_{j}(t_{i}) & \ \in \text{ker}\left( \Psi_{r,j}(t_{i},t_{0}) \right) \\
\Delta\bm{v}(t_{i}) & \ = \Psi_{v,j}(t_{i},t_{0})\Delta\bm{c}_{j}(t_{i}) \\
\bm{c}_{j}(t_{i-1}) & \ + \Delta\bm{c}_{j}(t_{i})  \in \mathcal{C}^{(N,j)}
\end{align}
\end{subequations}
The first constraint is that the position is not changed by the impulsive maneuver. The second constraint is true by definition of the fundamental solution matrix $\Psi_{j}$ in Cartesian coordinates. The third constraint states that regardless of the change induced in $\bm{c}_{j}$, the new monomial state must still be on the manifold $\mathcal{C}^{(N,j)}$. We can now rewrite the problem of Eq.~\eqref{TrajOpt1}:
\begin{equation}
\label{TrajOpt4}
\begin{split}
& \underset{\bm{c}_{j}(t_{i}) \ \forall \ i \in[1,K]}{\text{min}} \ \ \ \sum_{i=1}^{K}\left.\| \Psi_{v,j}(t_{i})\left( \bm{c}_{j}(t_{i}) - \bm{c}_{j}(t_{i-1}) \right) \|\right.^{\mathcal{P}} \\ \\
& \text{subject to} \ \ \ \begin{matrix} & \ \bm{c}_{j}(t_{0}) = \bm{c}_{j,\text{start}} \\ & \ \bm{c}_{j}(t_{K}) = \bm{c}_{j,\text{goal}}  \\ & \ \bm{c}_{j}(t_{i}) - \bm{c}_{j}(t_{i-1}) \in \text{ker}\left(\Psi_{r,j}(t_{i})\right) \\
& \ \bm{c}_{j}(t_{i}) \in \mathcal{C}^{(N,j)}, \ \ \forall i \geq 1
\end{matrix} 
\end{split}
\end{equation}
The problem is reduced to choosing reasonable discrete steps $\Delta\bm{c}_{j}(t_{i})$, $i = 1, 2, \ldots, K$, achieving the desired path from $\bm{c}_{j,\text{start}}$ at $t_{0}$ to $\bm{c}_{j,\text{goal}}$ at $t_{f}$, while minimizing the above cost. For $K\gg N$, the optimal solution will generally return many null jumps in the states, i.e. $\bm{c}_{j}(t_{q}) = \bm{c}_{j}(t_{q-1})$ for sub-optimal maneuver time $t_{q}$. Comparing Eq.~\eqref{TrajOpt1} and Eq.~\eqref{TrajOpt4}, it should be clear that the latter is an orverparameterized representation of the former. Note additionally that the dynamical constraints no longer appear in this new representation. The final state $\bm{c}_{j,\text{goal}}$ is obtained from $\bm{x}_{\text{goal}}(t_{f})$ via a single simple inversion of the map $\Psi_{j}(t_{f},0)$ (see e.g. \cite{Berz_ParticleBeamMaps}). In particular, the initial guess is given by inversion of the linear part of the map, which is just the STM:
\begin{equation}
\label{igGOAL}
\bm{c}_{1,\text{goal}}[0] = \Phi^{-1}(t_{f},0)\bm{x}_{\text{goal}}(t_{f})
\end{equation}
then subsequent nonlinearity corrections are implemented in a Newton-style framework: 
\begin{equation}
\label{DiffCorrGOAL}
\bm{c}_{1,\text{goal}}[p] = \bm{c}_{1,\text{goal}}[p-1] + \alpha\left(\Psi_{j}(t_{f})\left.\frac{\partial\bm{c}_{j}}{\partial\bm{c}_{1}}\right\vert_{\bm{c}_{1,\text{goal}}[k-1]}\right)^{-1}\left(\bm{x}_{\text{goal}}(t_{f}) - \Psi_{j}(t_{f})\bm{E}_{j}\left(\bm{c}_{1,\text{goal}}[k-1]\right)\right)
\end{equation}
where $p$ is the iteration number and $\alpha \leq 1$ scales the step for added stability. For our examples the simple choice $\alpha=1$ is sufficient, and the method is not numerically burdensome, typically executing in Python in $\sim0.01$ s for $N=6$, $j=3$. Upon convergence, we obtain $\bm{c}_{j,\text{goal}} = \bm{E}_{j}(\bm{c}_{1,\text{goal}})$.

The problem given by Eq.~\eqref{TrajOpt4} is intended as a simple example, but in general many functional constraints on the state $\bm{x}(t_{i})$ can be inherited as linear constraints on the $\bm{c}_{j}(t_{i})$. 
For this optimal path-planning problem, we write out a vector of (preliminary) decision variables as $\tilde{\bm{X}} = \left(\bm{c}_{j}(t_{1})^{\top}, \bm{c}_{j}(t_{2})^{\top}, \ldots, \bm{c}_{j}(t_{K})^{\top}\right)^{\top}$. The cost function of the problem given by Eq.~\eqref{TrajOpt4} can be shown to be quadratic in $\tilde{\bm{X}}$, and the first two constraints are clearly linear. The third constraint is also linear, written simply as $\Psi_{r,j}(t_{i})\Delta\bm{c}_{j}(t_{i}) = \bm{0}$. 
The only non-convex part of the problem given by Eq.~\eqref{TrajOpt4} is the final constraint that the decision variables lie on the manifold $\mathcal{C}^{(N,j)}$. For this non-convexity we propose a sequential convex programming problem.

\subsection{Sequential convex programming}
\subsubsection{Defining the SCP}
\begin{figure}[h!]
\centering
\includegraphics[]{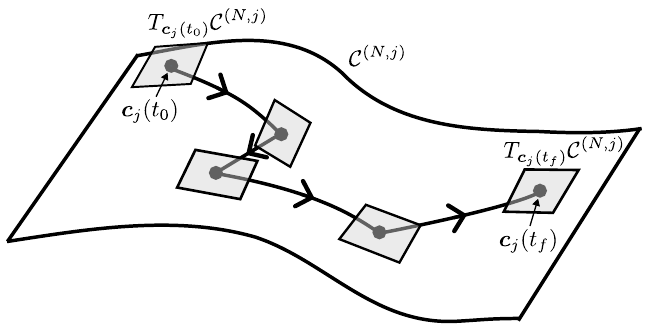}
\caption{Sequential Convexification via Monomial Coordinates on $\mathcal{C}^{(N,j)}$}
\label{fig:myfig2}
\end{figure}

Here we define the SCP approach used for generating solutions to numerical examples later in this work. Figure~\ref{fig:myfig2} conceptually depicts a trajectory, in terms of 5 distinct points in the monomial coordinates, resulting from 4 impulsive maneuvers. 
In the vicinity of any of the points in the trajectory shown in Fig.~\ref{fig:myfig2}, we can linearly approximate local variations in $\bm{c}_{j}$ to lie in the tangent space as below:
\begin{equation}
\label{scvx1}
\delta\bm{c}_{j}(t_{i}) \approx \left.\frac{\partial\bm{c}_{j}}{\partial\bm{c}_{1}}\right\vert_{\bm{c}_{j}(t_{i})}\delta\bm{c}_{1}(t_{i})  \ \in \ T_{\bm{c}_{j}(t_{i})}\mathcal{C}^{(N,j)}
\end{equation}
This expression is only true in the infinitesimal sense, but for a sufficiently smooth constraint surface it is a decent approximation for larger variations. In this case we know analytically the form of $\mathcal{C}^{(N,j)}$, the only non-convexity in our problem. The key to a sequential convex programming (SCP) implementation of the problem given by Eq.~\eqref{TrajOpt4} is to make use of the local tangent plane approximation of Eq.~\eqref{scvx1}, and also to choose our free variables for the convex problem as the $\delta\bm{c}_{1}(t_{i})$, with the nominal $\bm{c}_{j}(t_{i})$ chosen from a prior iteration (or, for iteration 1, from the initial guess of the trajectory). Thus $\tilde{\bm{X}} = (\delta\bm{c}_{1}(t_{1})^{\top}, \delta\bm{c}_{1}(t_{2})^{\top}, \ldots, \delta\bm{c}_{1}(t_{K})^{\top})^{\top}$, and transforming Eq.~\eqref{TrajOpt4}, the resulting convex sub-problem is 
given below:
\begin{equation}
\label{TrajOpt6}
\begin{split}
\allowdisplaybreaks
& \underset{\delta\bm{c}_{1}(t_{i}), \ \bm{s}_{i} \ \forall \ i \in[1, K], \ \bm{s}_{\text{end}}}{\text{min}} \ \ \ \begin{matrix} \left\|\Psi_{v,j}(t_{1})\left(\bm{c}_{j}^{\dagger}(t_{1}) + \left.\frac{\partial\bm{c}_{j}}{\partial\bm{c}_{1}}\right\vert_{1}\delta\bm{c}_{1}(t_{1})  - \bm{c}_{j}(t_{0}) \right) \right\|^{\mathcal{P}} + w\left(\sum_{i=1}^{K}\left\|\bm{s}_{i}\right\|^{2} + \|\bm{s}_{\text{end}}\|^{2}\right) \\ + \sum_{i=2}^{K} \left\|\Psi_{v,j}(t_{i})\left(\bm{c}_{j}^{\dagger}(t_{i}) + \left.\frac{\partial\bm{c}_{j}}{\partial\bm{c}_{1}}\right\vert_{i}\delta\bm{c}_{1}(t_{i})  - \bm{c}_{j}^{\dagger}(t_{i-1}) - \left.\frac{\partial\bm{c}_{j}}{\partial\bm{c}_{1}}\right\vert_{i-1}\delta\bm{c}_{1}(t_{i-1})  \right) \right\|^{\mathcal{P}} \end{matrix} \\ \\
& \text{subject to} \ \ \ \begin{matrix} 
\Psi_{j}(t_{k})\left(\bm{c}_{j}^{\dagger}(t_{K}) +  \left.\frac{\partial\bm{c}_{j}}{\partial\bm{c}_{1}}\right\vert_{K}\delta\bm{c}_{1}(t_{K})\right) + \bm{s}_{\text{end}} = \bm{x}_{\text{goal}}(t_{K})
\\ \Psi_{r,j}(t_{1})\left( \bm{c}_{j}^{\dagger}(t_{1}) + \left.\frac{\partial\bm{c}_{j}}{\partial\bm{c}_{1}}\right\vert_{1}\delta\bm{c}_{1}(t_{1}) - \bm{c}_{j}(t_{0}) \right) + \bm{s}_{1} = \bm{0}  \\ \Psi_{r,j}(t_{i+1})\left( \bm{c}_{j}^{\dagger}(t_{i+1}) - \bm{c}_{j}^{\dagger}(t_{i}) + \left.\frac{\partial\bm{c}_{j}}{\partial\bm{c}_{1}}\right\vert_{i+1}\delta\bm{c}_{1}(t_{i+1}) - \left.\frac{\partial\bm{c}_{j}}{\partial\bm{c}_{1}}\right\vert_{i}\delta\bm{c}_{1}(t_{i}) \right)+ \bm{s}_{i+1} = \bm{0}, i = 1\text{:}K-1 
\end{matrix} 
\end{split}
\end{equation}
Note the careful parsing of expressions involving $\bm{c}_{j}(t_{0})$, which is fixed in this example and is not part of the decision variables. The $( \ )^{\dagger}$ denotes terms from the solution to the prior iteration, about which the current iteration is expanded. Furthermore $\left.\frac{\partial\bm{c}_{j}}{\partial\bm{c}_{1}}\right\vert_{i}$ is a shorthand for $\left.\frac{\partial\bm{c}_{j}}{\partial\bm{c}_{1}}\right\vert_{\bm{c}_{j}^{\dagger}(t_{i})}$. The enumerated slack variables $\bm{s}_{i}$ are introduced to prevent artificial infeasibility of the convex sub-problem, with a scalar weight $w>0$ to specify the degree of penalization of slack terms (the converged solution must satisfy $\bm{s}_{i} = \bm{0}$). These slack variables can be physically interpreted as the positional defect constraints in the trajectory. There is also an additional slack variable $\bm{s}_{\text{end}}$ related to the satisfaction of the first listed linear constraint, which is a constraint on the end state, and is similarly penalized. For performance reasons, it is best for the problem states to be rendered non-dimensional, or for the last $N/2$ components of $\bm{s}_{\text{end}}$ to be rescaled to have the same positional ``units" as the $\bm{s}_{i}$ slack variables. In this manner an SCP iteration will not be biased to over/under penalize any components $\bm{s}_{\text{end}}$ in comparison to the other $\bm{s}_{i}$. In our experience with this formulation, proper numerical scaling of the free variables and constraints is necessary to ensure good performance when using commercial or open-source solvers.

Between iterations, it is necessary to project the solution to the convex sub-problem (which exists in the union of the tangent planes of points $\bm{c}_{j}^{\dagger}(t_{i}) \ \forall \ i$) back onto the manifold. The most obvious (and computationally easiest) way is to compute for each update $+\delta\bm{c}_{1}(t_{i})$ the new $\bm{c}_{j}(t_{i})$ given by $\bm{c}_{1}^{\dagger}(t_{i}) + \delta\bm{c}_{1}(t_{i})$, e.g. using Eq.~\eqref{Ej}.

Because the tangent-plane approximations of variations $\delta\bm{c}_{j}$ are only locally valid, we must additionally impose some kind of trust region constraint preventing the sub-problem from obtaining variations $\delta\bm{c}_{1}(t_{i})$ that are too large. In this work we impose a norm constraint on $\tilde{\bm{X}}$ to do this:
\begin{equation}
\label{NormConstr1}
\|\tilde{\bm{X}}\| \leq d
\end{equation}
We can set a fixed trust-region radius $d$, or alternatively we can update this via a proper trust-region update method such as the one outlined in \cite{HofmannDART2022scitech}. 

\subsubsection{Parsing the convex sub-problem for $\mathcal{P}=2$}
The sub-problem defined by Eq.~\eqref{TrajOpt6} can be resolved for $\mathcal{P}=2$ in the form below:
\begin{equation}
\label{QP1}
\begin{split}
& \underset{\bm{X}}{\text{min}} \ \  \bm{X}^{\top}P\bm{X} + \bm{q}^{\top}\bm{X} \\ \\
& \text{subject to} \ \ \begin{matrix}  A\bm{X} = \bm{b} \\ \| M \bm{X} \| \leq d  \end{matrix} 
\end{split}
\end{equation}
This parsing is provided for implementational convenience. Note however that our numerical results use CVXPY to automatically parse the original sub-problem as a second-order cone program (SOCP), then solve via ECOS. The form of the above convex problem changes depending on the choice of norm for the final constraint. Choosing a 2-norm, it becomes a quadratically constrained quadratic program (QCQP). Choosing instead the infinity-norm, it can be written easily as a quadratic problem (QP) by expanding out the resulting absolute value constraints on each component of $\bm{X}$, yielding  $M_{2}\bm{X} \leq \bm{d}_{2}$ (see e.g. \cite{Boyd_Vandenberghe_2004}). This is lower than QCQP and SOCP on the hierarchy of convex optimization problems, and should solve more quickly. 

For the above formulation, we first augment the preliminary decision variables $\tilde{\bm{X}}$ with slack variables $\bm{S} = (\bm{s}_{1}^{\top}, \bm{s}_{2}^{\top}, \ldots \bm{s}_{K}^{\top}, \bm{s}_{\text{end}}^{\top})^{\top}$ to form $\bm{X} = (\tilde{\bm{X}}^{\top}, \bm{S}^{\top})^{\top}$, which is of length $\left(\frac{3}{2}K + 1\right) N$ for assumed even $N$. Because we only want to constrain $\tilde{\bm{X}}$, the form of $M$ is simple:
\begin{equation}
\label{QP2}
M = \begin{bmatrix} I_{NK\times NK} \ \ 0_{NK \times N\left(\frac{K}{2} + 1\right)} \end{bmatrix}
\end{equation} 
Because all constraints in Eq.~\eqref{TrajOpt6} are linear and fairly simple, it is easy to construct $A$ and $\bm{b}$. The form of $P$ and $\bm{q}$ however requires some algebra and will be provided. 

While it does not influence the convex sub-problem, we must compute the part of the cost $J$ in Eq.~\eqref{TrajOpt6} that is not a function of the decision variables of a given iteration. This term is in fact the total cost predicted from the prior iteration, so we denote it as $J^{\dagger}$, reusing our dagger notation:
\begin{subequations}
\label{MakeJ}
\begin{align}
J^{\dagger} = & \ \sum_{i=i}^{K}\bm{c}_{j}^{\dagger\top}(t_{i})O_{i}\bm{c}_{j}^{\dagger}(t_{i}) - 2\bm{c}_{j}^{\dagger\top}(t_{i})O_{i}\bm{c}_{j}^{\dagger}(t_{i-1}) + \bm{c}_{j}^{\dagger\top}(t_{i-1})O_{i}\bm{c}_{j}^{\dagger}(t_{i-1}) \\
O_{i} = & \ \Psi_{v,j}^{\top}(t_{i})\Psi_{v,j}(t_{i})
\end{align}
\end{subequations}
Thus the solution to Eq.~\eqref{QP1} provides the minimal admissible $\delta J =  \bm{X}^{\top}P\bm{X} + \bm{q}^{\top}\bm{X}$ yielding total cost $J^{(q)} \approx J^{(q-1)} + \delta J^{(q)}$ for iteration $q$. In reality, after iteration $q$, the nonlinear projection of the linear step yields the true $J^{(q)}$ after a procedure of 1) nonlinearly updating all $\bm{c}_{j}^{\dagger}(t_{i})$, 2) setting all $\delta\bm{c}_{1}(t_{i}) = \bm{0}$ in Eq.~\eqref{TrajOpt6}, then 3) using that to compute the true values of all slack variables, which are simply the true defects in the constraint equations after the nonlinear projection. From this info, we can compute the actual increment in cost $\delta J$. We can use the disparity between the predicted (linear) $\delta J$ and the actual (nonlinear) value to determine how nonlinear an SCP iteration was, with greater nonlinearity yielding greater disparity. The two quantities are given below:
\begin{subequations}
\label{dJEq}
\begin{align}
\delta J_{\text{predict}}^{(q)} = & \ \bm{X}^{\top}P\bm{X} + \bm{q}^{\top}\bm{X} \\
\delta J_{\text{actual}}^{(q)} = & \ J^{(q)} - J^{(q-1)}
\end{align}
\end{subequations}
where $P$, $\bm{q}$ are computed using information from iteration $q-1$.

We now define the matrix $P$ and vector $\bm{q}$ by sequential accounting instead of explicitly writing the large matrices. In this manner, they are initialized as zeros with the correct size, and consideration of each successive cost term in Eq.~\eqref{TrajOpt6}  adds a corresponding part to the matrices whose form is to be defined.
\begin{subequations}
\label{MakeMats0}
\allowdisplaybreaks
\begin{align}
P = & \ P_{0} + \sum \delta P, \ \ \ P_{0} = 0_{L\times L} \\
\bm{q} = & \ \bm{q}_{0} + \sum \delta \bm{q}, \ \ \ \bm{q}_{0} = \bm{0}_{L\times 1}  \\
L = & \ \left(\frac{3}{2}K + 1\right) N
\end{align}
\end{subequations}
Starting with $P$, the first ($t_{1}$), second (slack-associated), and third (indexed) cost terms of Eq.~\eqref{TrajOpt6}  produce the following contributions to $P$:
\begin{subequations}
\label{MakeP1}
\begin{align}
& \delta P_{1} = \left[ \begin{array}{c; {2pt/2pt}c} C_{1}^{\top}O_{1}C_{1} & 0_{N\times(L-N)} \\ \hdashline[2pt/2pt] 0_{(L-N)\times N} & 0_{(L-N)\times(L-N)} \end{array}\right] \\
& C_{1} = \left.\frac{\partial\bm{c}_{j}}{\partial\bm{c}_{1}}\right\vert_{\bm{c}_{j}^{\dagger}(t_{1})}, \ \ \ O_{1} = \Psi_{v,j}^{\top}(t_{1})\Psi_{v,j}(t_{1})
\end{align}
\end{subequations}
\begin{subequations}
\label{MakeP2}
\begin{align}
& \delta P_{2} =  \left[ \begin{array}{c; {2pt/2pt}c}  0_{NK\times NK} & 0_{NK \times L_{S}} \\ \hdashline[2pt/2pt] 0_{L_{S}\times NK} & w I_{L_{S}\times L_{S}} \end{array}\right]  \\
& L_{S} = L - KN = N\left(\frac{K}{2} + 1\right)
\end{align}
\end{subequations}
\begin{subequations}
\label{MakeP3}
\begin{align}
& \delta P_{3,i} = \left[ \begin{array}{c; {2pt/2pt}c; {2pt/2pt}c; {2pt/2pt}c} 0_{\mathcal{A}\times \mathcal{A}} & 0_{\mathcal{A} \times N} & 0_{\mathcal{A}\times N} & 0_{\mathcal{A}\times\mathcal{B}} \\ \hdashline[2pt/2pt] \\[-1em] 0_{N\times \mathcal{A}} & C_{i-1}^{\top}O_{i}C_{i-1} & -C_{i-1}^{\top}O_{i}C_{i} & 0_{N\times\mathcal{B}} \\ \hdashline[2pt/2pt] \\[-1em] 0_{N\times \mathcal{A}} & -C_{i}^{\top}O_{i}C_{i-1} & C_{i}^{\top}O_{i}C_{i} & 0_{N\times\mathcal{B}} \\ \hdashline[2pt/2pt] 0_{\mathcal{B}\times\mathcal{A}} & 0_{\mathcal{B}\times N} & 0_{\mathcal{B}\times N} & 0_{\mathcal{B}\times\mathcal{B}} \end{array}\right] \\
& C_{i} = \left.\frac{\partial\bm{c}_{j}}{\partial\bm{c}_{1}}\right\vert_{\bm{c}_{j}^{\dagger}(t_{i})}, \ \ \ O_{i} = \Psi_{v,j}^{\top}(t_{i})\Psi_{v,j}(t_{i}) \\
& \mathcal{A} = N(i-2), \ \ \ \mathcal{B} = L - Ni, \ \ \ i = 2\text{:}K 
\end{align}
\end{subequations}
where the $\delta P_{3,i}$ is added for every applicable index $i$. Now $\bm{q}$ is defined similarly for the first ($t_{1}$) and third (indexed) cost terms:
\begin{equation}
\label{dq1}
\delta\bm{q}_{1} =  \left[\begin{array}{c; {2pt/2pt} c} 2\left(\bm{c}_{j}^{\dagger}(t_{1}) - \bm{c}_{j}(t_{0})\right)^{\top}O_{1}C_{1} & \bm{0}_{1\times L-N} \end{array}\right]^{\top} 
\end{equation}
\begin{equation}
\label{dq2}
\delta\bm{q}_{2,i} = \left[\begin{array}{c; {2pt/2pt} c; {2pt/2pt} c; {2pt/2pt} c} \bm{0}_{1\times\mathcal{A}} & -2\left(\bm{c}_{j}^{\dagger}(t_{i}) - \bm{c}_{j}^{\dagger}(t_{i-1})\right)^{\top}O_{i}C_{i-1} & 2\left(\bm{c}_{j}^{\dagger}(t_{i}) - \bm{c}_{j}^{\dagger}(t_{i-1})\right)^{\top}O_{i}C_{i} &  \bm{0}_{1\times\mathcal{B}} \end{array} \right]^{\top}, \ \ \ i = 2\text{:}K 
\end{equation}

Beyond the nonlinear projection, the only other nonlinear calculations needed for the setup of Eq.~\eqref{TrajOpt6} between iterations are 1) the analytic re-computation of all $\left.\frac{\partial\bm{c}_{j}}{\partial\bm{c}_{1}}\right\vert_{i}$ Jacobians and 2) the computation of any needed norms. 
This simple implementation liberates the details of dynamics from the problem representation. The converged solution is guaranteed to retain the full accuracy of the $j$\textsuperscript{th}-order fundamental solution expansion because it serves as our transcription scheme. Furthermore, given a measure of the domain of validity $\mathcal{C}_{\varepsilon}$, it is straightforward to enforce trust region radii $d_{i}$ in the sub-problem such that at all times each $\|\bm{c}_{1}^{\dagger}(t_{i}) + \delta\bm{c}_{1}(t_{i})\|$ for each trajectory node never extends outside the domain of validity, ensuring robustness of the guidance scheme. 
\subsection{Other solution strategies}
The SCP scheme that we define for solving the path-planning problem on the manifold is provided as a proof-of-concept. Here, for completeness, we provide a very brief summary of other strategies that can be adapted to solve the problems defined in Eqs.~\eqref{TrajOpt4} and \eqref{TrajOpt6}, exploiting the structure of $\mathcal{C}^{(N,j)}$. First, we have already mentioned that $\mathcal{C}^{(N,j)}$ is an embedded submanifold of the Euclidean space $\mathbb{R}^{K_{j}}$. Adopting the induced metric from $\mathbb{R}^{K_{j}}$, $\mathcal{C}^{(N,j)}$ becomes a Riemannian submanifold \citep{boumal2023intromanifolds}. Then, solving for a sequence of points $\left(\bm{c}_{j}(t_{1}), \bm{c}_{j}(t_{2}), \ldots, \bm{c}_{j}(t_{K})\right)$ can be reformulated as solving for a single point on the $KN$-dimensional \textit{product manifold} $\mathcal{C}^{(N,j)} \times \mathcal{C}^{(N,j)} \times \ldots \times \mathcal{C}^{(N,j)}$ subject to the reformulated kinematic constraints, which will act as an equality constraint on the product state. \cite{RiemannianIP} discusses the development of an interior point method for use on Riemannian manifolds, which could in principle solve the problem posed by Eq.~\eqref{TrajOpt4}, including accommodation of equality and inequality state constraints. Furthermore, absorbing the equality constraint (see e.g. \cite{Boyd_Vandenberghe_2004}), the resulting unconstrained problem could in principle be solved by the Riemannian Trust Regions (RTR) algorithm of \cite{boumal2023intromanifolds}, developed specifically for Riemannian sub-manifolds. Both of these methods have the benefits of proof of convergence. 


\section{Example Problem: Spacecraft Rendezvous and Station-Keeping}
In this section, we showcase the potential use of the developments in this paper for an example application from space vehicle guidance, to establish both practical and theoretical interest in this overparameterized monomial formulation.
\subsection{Nonlinear spacecraft rendezvous in LEO -- Two-stage guidance}
Before showing results with an SCP implementation, we demonstrate a two-stage guidance scheme (linear predictor, nonlinear corrector) for rapid generation of trajectories that are both fuel-efficient and dynamically feasible. For weakly nonlinear cases, this scheme extends the reach of linearized guidance solutions for a small delta-V penalty. For more strongly nonlinear cases, it provides a convenient feasible (but sub-optimal) initial guess for the SCP solver. The methods in this paper can be applied in to rendezvous and station-keeping in practically any dynamical environment, but we consider the classical rendezvous equations with a target in a simple circular Earth orbit due to its historical and modern significance for rendezvous. 
Using the framework in this paper, we demonstrate the following convenient two-stage guidance strategy for spacecraft rendezvous:
\begin{itemize}
\item [] \textbf{Stage 1:} Compute a low-fidelity solution by solving the linearized problem (i.e. $j=1$). This can be done using several state-of-the-art convex optimization-based approaches.
\item [] \textbf{Stage 2:} Rapid nonlinear correction, to accommodate for nonlinearity effects. This is done by leveraging computation of the high-fidelity fundamental solution matrix $\Psi_{j}(t)$ and the analytic expression for the monomial vector $\bm{c}_{j}$ at desired order of nonlinearity $j$. 
\end{itemize}
The nonlinear fundamental solutions data necessary for Stage 2 can be computed in advance and stored on a spacecraft's guidance computer. With this approach, implementation of Stage 2 is reduced to the use of pre-programmed analytic functions, data lookup, and manageable linear algebra. In settings where many trajectories need to be computed and characterized, pre-computation of the nonlinear fundamental solution data is substantially computationally liberating. 

\subsubsection{Stage 1}
First, a suitable linearized rendezvous guidance scheme is defined. The nonlinear dynamics under consideration are those of relative motion in the vicinity of a simple circular Keplerian orbit. This problem is parameterized by the nonlinear equations below for LVLH relative position $\bm{r} = (x,y,z)^{\top}$ and velocity $\bm{v} = (\dot{x},\dot{y},\dot{z})$ \citep{schaub}: 
\begin{subequations}
\label{SchaubFFNL}
\begin{align}
\ddot{x} - 2n\dot{y} - xn^{2} - \frac{\mu}{R^{2}} = & \ -\frac{\mu}{r^{3}}(R+x) \\
\ddot{y} + 2n\dot{x} - yn^{2} = & \ -\frac{\mu}{r^{3}}y \\
\ddot{z} = & \ -\frac{\mu}{r^{3}}z
\end{align}
\end{subequations}
where $R$ is the circular target orbit radius, $n$ is its orbital mean motion, $r = \sqrt{(R+x)^{2} + y^{2} + z^{2}}$ is the chaser's instantaneous orbit radius, and velocities are measured with respect to the target-centered LVLH frame. Note for this example that $\bm{X}_{r}(t)$ is the orbital state of the target spacecraft, but the formulation of Eq.~\eqref{SchaubFFNL} describes only the dynamics of the relative states.

In the vicinity of the target spacecraft, the nonlinear Eq.~\eqref{SchaubFFNL} linearizes to the popular Clohessy-Wiltshire (CW) model \citep{clohessy1}:
\begin{subequations}
\label{CW1}
\begin{align}
\ddot{x} = & \ 2n\dot{y} + 3n^{2}x\\
\ddot{y} = & \ -2n\dot{x} \\
\ddot{z} = & \ -n^{2}z
\end{align}
\end{subequations}
We have explored solving this linear problem in two ways. First, reusing a method from \cite{BurnettJGCD2022} which indirectly solves for the optimal burn times and directions, then linearly solves for the resultant magnitudes along each direction (see also \cite{Guffanti_ICs} and references therein). Second, using a slightly slower method that is more stable and more clear in the context of this work, which is defined below as the linearized equivalent of the problem given by Eq.~\eqref{TrajOpt4}:
\begin{equation}
\label{Stage1LINEAR}
\begin{split}
& \underset{\bm{c}_{1}(t_{i}) \ \forall i \in [1,K]}{\text{min}} \ \ \sum_{i=1}^{K}\left.\| \Phi_{v}(t_{i})\left( \bm{c}_{1}(t_{i}) - \bm{c}_{1}(t_{i-1}) \right) \|\right.^{\mathcal{P}} \\ \\
& \text{subject to} \ \ \begin{matrix} \bm{c}_{1}(t_{0}) = \bm{c}_{1,\text{start}} \\ \bm{c}_{1}(t_{K}) = \bm{c}_{1,\text{goal}}  \\ \bm{c}_{1}(t_{i}) - \bm{c}_{1}(t_{i-1}) \in \text{ker}\left(\Phi_{r}(t_{i})\right) 
\end{matrix} 
\end{split}
\end{equation}
where $\Phi_{r} = \Psi_{r,1}$ and $\Phi_{v} = \Psi_{v,1}$ denote the top and bottom halves of the state transition matrix solution to the linearized dynamics in Eq.~\eqref{CW1}, and they can also be thought of as just the $j=1$ part of the map $\Psi_{j}$ computed from Eq.~\eqref{SchaubFFNL}. This linearized problem is already convex, for which we define the free variable vector $\bm{X} = (\bm{c}_{1}^{\top}(t_{1}), \bm{c}_{1}^{\top}(t_{2}), \ldots, \bm{c}_{1}^{\top}(t_{K}))^{\top}$. The solver reveals the $k$ optimal discretized burn times $t_{b_{i}}$ as those for which non-null jumps in the states are computed, $\bm{c}_{1}(t_{b_{i}}) - \bm{c}_{1}(t_{b_{i-1}}) \neq \bm{0}$. By contrast, for all times $t_{q}$ that would yield sub-optimal burns, the solution to Eq.~\eqref{Stage1LINEAR} will compute $\bm{c}_{1}(t_{q}) - \bm{c}_{1}(t_{q-1}) = \bm{0}$ to the precision of the convex solver used.

The solution of Stage 1 consists of the following:
\begin{enumerate}
\item A sequence of (linear) monomial states $\bm{c}_{1}(\tau_{i})$ at $k+1$ times $\tau_{i} \in \mathcal{T}$. 
\item A list of critical times $\mathcal{T} = [t_{0}, t_{b_{1}}, t_{b_{2}}, \ldots, t_{b_{k}}]$, where $\tau_{0} = t_{0}$ is the initial time and $\tau_{i} = t_{b_{i}}$ is the $i$\textsuperscript{th} burn time. With this notation, $\bm{c}_{1}(t_{0})$ represents the epoch state before the control action of the first burn at $t_{b_{i}} \geq t_{0}$ is taken. 
\item A table of $k$ associated nonzero delta-V vectors. 
\end{enumerate}
The outputs of any other guidance algorithm based on linearization and the impulsive maneuver assumptions (see e.g. \cite{berning2024BO} for a recent example of linearized guidance incorporating the practical constraint of passive safety) can also be put into this form. Note the use of lower-case ``k" instead of ``K" from Section III because in this context each increment of the index is by definition a (nonzero) maneuver, selected from a much larger set of candidate burn times, thus $k \ll K$.

\subsubsection{Stage 2}
This scheme is a rapid high-fidelity correction to the solution of Stage 1 which preserves the burn times $t_{b_{i}} \in \mathcal{T}_{b}$ and controllable maneuver positions $\bm{r}(t_{b_{i}})$, while altering the magnitude and direction of the delta-Vs $\Delta\bm{v}(t_{b_{i}})$. The result is a guidance solution compatible with the dynamics parameterized by a given nonlinear fundamental solution representation $\Psi_{j}(t)$ and associated overparameterized monomial set $\bm{c}_{j}$. The solution is trustworthy if the nonlinear fundamental solution parameterization $\Psi_{j}(t)$ is accurate for the flight environment and domain under consideration.

For this scheme simple Newton-style correction is sufficient. We define a free variable vector $\bm{X}$ in terms of the (first-order) monomials $\bm{c}_{1}(t_{b_{i}})$, and a constraint vector $\bm{F}(\bm{X})$, which is driven to zero as the satisfactory $\bm{X}^{*}$ is obtained. This function is linear in the higher-order monomial representation $\bm{c}_{j}(t_{b_{i}})$ and nonlinear in the $\bm{c}_{1}(t_{b_{i}})$. We demonstrate now the form of $\bm{X}$ and $\bm{F}(\bm{X})$:
 \begin{equation}
 \label{2SGpost1}
 \bm{X} = \left( \bm{c}_{1}(t_{b_{1}})^{\top}, \bm{c}_{1}(t_{b_{2}})^{\top}, \bm{c}_{1}(t_{b_{3}})^{\top}, \ldots, \bm{c}_{1}(t_{b_{k}})^{\top} \right)^{\top}
 \end{equation}
 
 \begin{equation}
 \label{2SGpost2}
 \bm{F}(\bm{X}) = \begin{pmatrix} \Psi_{r,j}(\tau_{2})\bm{c}_{j}(\tau_{2}) - \bm{r}_{G}(\tau_{2}) \\ \Psi_{r,j}(\tau_{3})\bm{c}_{j}(\tau_{3}) - \bm{r}_{G}(\tau_{3}) \\ \vdots \\ \Psi_{r,j}(\tau_{k})\bm{c}_{j}(\tau_{k}) - \bm{r}_{G}(\tau_{k}) \\ \hdashline[2pt/2pt] \\[-2ex] \Psi_{r,j}(\tau_{1})\bm{c}_{j}(\tau_{1}) - \Psi_{r,j}(\tau_{1})\bm{c}_{j}(\tau_{0}) \\ \Psi_{r,j}(\tau_{2})\left( \bm{c}_{j}(\tau_{2}) - \bm{c}_{j}(\tau_{1})\right) \\ \Psi_{r,j}(\tau_{3})\left( \bm{c}_{j}(\tau_{3}) - \bm{c}_{j}(\tau_{2})\right) \\ \vdots \\ \Psi_{r,j}(\tau_{k})\left( \bm{c}_{j}(\tau_{k}) - \bm{c}_{j}(\tau_{k-1})\right) \\ \hdashline[2pt/2pt] \\[-2ex] \Psi_{v,j}(\tau_{k})\bm{c}_{j}(\tau_{k}) - \bm{v}_{G}(\tau_{k}) \end{pmatrix}
 \end{equation}
  
where the nonlinear mapping from $\bm{c}_{1}(\tau_{i})$ (and hence from $\bm{X}$) to $\bm{c}_{j}(\tau_{i})$ is analytic and straightforward -- recall the discussion in Section II.B. The first set of constraints (above the dashed line) are enforcing that all controllable maneuver locations match the corresponding maneuver location from the guidance solution of Stage 1, $\bm{r}_{G}(\tau_{i})$. Because the first maneuver location is not controllable, it is omitted from these constraints. Note that if $\tau_{k} < t_{f}$, $\bm{x}_{G}(\tau_{k})$ (the state immediately after at the $k$\textsuperscript{th} and final burn) is related to the goal final condition $\bm{x}(t_{f})$ simply by the natural dynamics mapping $\bm{x}(\tau_{k}) = \bm{\psi}^{-1}(\bm{x}(t_{f}), t_{f}, \tau_{k})$. The second set of constraints, between dashed lines, enforce the kinematic constraint that changes in the monomial states at burn nodes (implicitly the result of impulsive maneuvers) should not instantaneously change the position. Note that the very first one is written differently to emphasize that $\bm{c}_{j}(\tau_{0})$ is uncontrollable and thus does not appear in the free variables $\bm{X}$. Lastly, the final listed constraint in $\bm{F}(\bm{X})$ is enforcing that the final goal velocity, immediately after the $k$\textsuperscript{th} and final burn, is achieved.

The constraint Jacobian, $G(\bm{X}) = \frac{\partial}{\partial\bm{X}}\left( \bm{F}(\bm{X}) \right)$, is computed analytically, and its derivation from Eq.~\eqref{2SGpost2} is straightforward:

\begin{equation}
 \label{2SGpost3}
 G(\bm{X}) = \begin{bmatrix} 0_{\frac{N}{2}\times N} & \Psi_{r,j}(\tau_{2})C_{j,2} & 0_{\frac{N}{2}\times N} & 0_{\frac{N}{2}\times N} & \ldots & 0_{\frac{N}{2}\times N} & 0_{\frac{N}{2}\times N} \\ 0_{\frac{N}{2}\times N} & 0_{\frac{N}{2}\times N} & \Psi_{r,j}(\tau_{3})C_{j,3} & 0_{\frac{N}{2}\times N} & \ldots & 0_{\frac{N}{2}\times N} & 0_{\frac{N}{2}\times N} \\ & & \vdots \\ 0_{\frac{N}{2}\times N} & 0_{\frac{N}{2}\times N} & 0_{\frac{N}{2}\times N} & 0_{\frac{N}{2}\times N} & \ldots & 0_{\frac{N}{2}\times N} & \Psi_{r,j}(\tau_{k})C_{j,k} \\[+1ex] \hdashline[2pt/2pt] \\[-2.5ex]  \Psi_{r,j}(\tau_{1})C_{j,1} & 0_{\frac{N}{2}\times N} & 0_{\frac{N}{2}\times N} & 0_{\frac{N}{2}\times  N} & \ldots & 0_{\frac{N}{2}\times N} & 0_{\frac{N}{2}\times N} \\ -\Psi_{r,j}(\tau_{2})C_{j,1} & \Psi_{r,j}(\tau_{2})C_{j,2} & 0_{\frac{N}{2}\times N} & 0_{\frac{N}{2}\times  N} & \ldots & 0_{\frac{N}{2}\times N} & 0_{\frac{N}{2}\times N} \\ 0_{\frac{N}{2}\times N} & -\Psi_{r,j}(\tau_{3})C_{j,2} & \Psi_{r,j}(\tau_{3})C_{j,3} & 0_{\frac{N}{2}\times N} & \ldots & 0_{\frac{N}{2}\times N} & 0_{\frac{N}{2}\times N} \\ & & \vdots \\ 0_{\frac{N}{2}\times N} & 0_{\frac{N}{2}\times N} & 0_{\frac{N}{2}\times N} & 0_{\frac{N}{2}\times N} & \ldots & -\Psi_{r,j}(\tau_{k})C_{j,k-1} & \Psi_{r,j}(\tau_{k})C_{j,k} \\ 0_{\frac{N}{2}\times N} & 0_{\frac{N}{2}\times N} & 0_{\frac{N}{2}\times N} & 0_{\frac{N}{2}\times N} & \ldots & 0_{\frac{N}{2}\times N} & \Psi_{v,j}(\tau_{k})C_{j,k} \end{bmatrix}
\end{equation}

where $C_{j,i} = \left.\frac{\partial\bm{c}_{j}}{\partial\bm{c}_{1}}\right\vert_{\bm{c}_{1}(\tau_{i})}$ which can be computed analytically by exploiting the known monomial structure of  the $\bm{c}_{j}$. This $kN \times kN$ matrix is full rank for our problem of interest. Thus the Newton step can be assessed as below:
\begin{equation}
\label{2SGpost4}
\delta\bm{X} = -\gamma G^{-1}\bm{F}(\bm{X}), \ \ 0 < \gamma \leq 1
\end{equation}
The proper $\gamma$ can, if needed, be chosen by evaluating the error norm of $\bm{F}(\bm{X} + \delta\bm{X})$ for a pre-determined sequence of $\gamma$, i.e. in a line-search, because of the extremely low overhead in computing $\bm{F}(\bm{X})$. For our examples, the simple choice $\gamma = 1$ was always sufficient. This linear problem is initialized by constructing $\bm{X}^{(0)}$ (i.e. the $0$\textsuperscript{th} iteration) directly from the $\bm{c}_{1}(\tau_{i})$ outputs from Stage 1 of guidance. 
Because there is no need for any complicated numerics in the analytic computations of Eqs.~\eqref{2SGpost3} and \eqref{2SGpost4}, Stage 2 is extremely rapid. In our experience, the total time taken by all needed successive corrections of Stage 2 are much faster than the small SOCP solved by stage 1.

We now apply the aforementioned procedure to an example short-range low-Earth orbit (LEO) rendezvous and proximity operations scenario. For this example, a linearized rendezvous trajectory is generated which requires correction for the nonlinear dynamics of relative motion in low-Earth orbit. The problem details are provided in Table~\ref{table:LEO1}. Also included in that table are the runtime details for generation of the nonlinear fundamental solutions in $\Psi_{j}$ via a differential algebra scheme. These are computed once, and reused for the first two examples in this paper. Note that instead using numerically propagated STTs, the compute times for each order were much slower, $\sim$ 129s (3), 6.69s (2), and 0.747s (1). Also, for this particular problem there exist analytic nonlinear expansions in literature. In particular, exact nonlinear expansions suitable for computing $\Psi_{j}(t)$ up to order $j=3$ can be found (in normalized form) in \cite{Butcher2016_AAS}. See also \cite{ButcherBurnett2017} and \cite{Willis2ndOrder} and references therein. Such analytic solutions unburden onboard algorithms from any onboard numerical integration to render the nonlinear fundamental solutions. More details can be found in the Appendix C.

\begin{table}[h!]
\centering
\caption{Parameters for Rendezvous Example 1}
\label{table:LEO1}
\begin{tabular}{ll} 
\hline\hline
\textbf{Category} & \textbf{Values} \\ 
\hline
\textbf{Target spacecraft} & LEO orbit, $a = 6378$ km, $e = 0$, $T \approx 1.4$ hrs \\ 
\textbf{Control interval} & $t_{0} = 0.1T$, $t_{f} = 2.3T$, discretized into 220 times \\ 
\textbf{Relative state at $t = 0$} & $\bm{x}(0) = [-1266.6, -12000, 1000, 0, 2.9748, 0]$ (units: m, m/s) \\ 
\textbf{Relative state at $t_{f}$} & $\bm{x}(t_{f}) = [-589.6, 383.2, -1825.9, 2.3747, 1.4617, -1.3499]$ \\ 
\hline
\textbf{STT info} & \\
Duration & $2.3T$ \\
Discretization & 230 times \\ 
Size & (order $j=3$): 474 kb \\ 
Compute time (order) & 0.795s (3), 0.349s (2), 0.175s (1) \\
\hline
\textbf{Linear guidance} & \\
Burn indices & $[0, 93, 142, 201, 219]$ \\ 
Delta-V & 2.336 m/s \\ 
Runtime & 0.108s \\ 
\hline
\textbf{Stage 2 guidance} & \\
Delta-V & 2.353 m/s \\ 
Runtime & 0.0119s \\ 
\hline\hline
\end{tabular}
\end{table}

\begin{figure}[h!]
\centering
\includegraphics[]{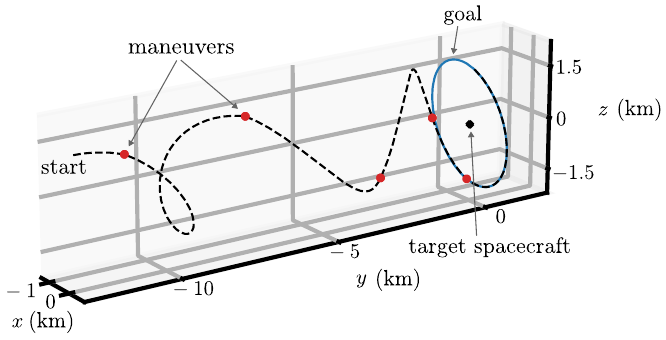}
\caption{Nominal linear guidance trajectory, linearized dynamics (rendezvous ex. 1)}
\label{fig:LEOGuidance_Linear}
\end{figure}

Figure~\ref{fig:LEOGuidance_Linear} gives the nominal guidance solution produced by Stage 1, propagated with the linearized dynamics assumed therein (in this case, the CW dynamics). This is generated on a 2023 MacBook Pro (M2 chip) in Python 3, using CVXPY \citep{diamond2016cvxpy,agrawal2018rewriting} with ECOS \citep{ECOS_2013} to solve the SOCP in $\sim$0.1s. After a sequence of 5 burns, the goal state is achieved: entry onto a tilted target-centered relative orbit. This linearized guidance solution requires a maneuver sequence with a total delta-V of 2.336 m/s.
In subsequent figures the nominal linear guidance solution appears as a gray curve. 

\begin{figure}[h!]
\centering
\includegraphics[]{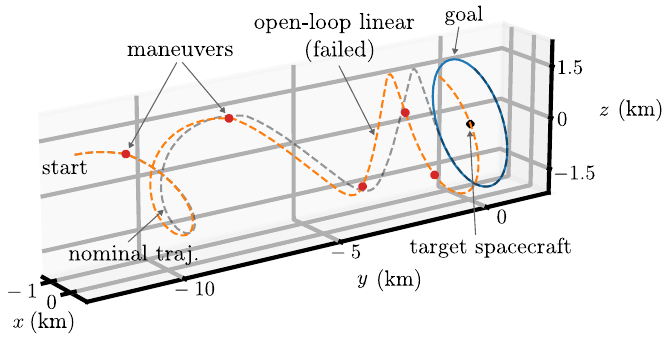}
\caption{Executed open-loop linear guidance, nonlinear dynamics (rendezvous ex. 1)}
\label{fig:LEOGuidance_Plot2}
\end{figure}

To produce Fig.~\ref{fig:LEOGuidance_Plot2}, the original Stage 1 guidance solution (nominal delta-Vs applied at pre-planned times) is followed in an open-loop fashion subject to the nonlinear dynamics -- in other words, the initial condition is integrated and the nominal impulsive guidance strategy is followed at the pre-computed maneuver times by stopping integration, adding the impulse, and resuming integration to the next impulse. This yields the orange curve. With no allowance for corrective maneuvers, the goal state is clearly not achieved. This illustrates a well-known limitation of using linearized dynamics for rendezvous guidance: the nominal trajectories predicted by linearization will have to be corrected for accuracy in a high-fidelity model. There are various means of doing this, such as refinement via sequential convexification, post-processing with sequential Lambert targeting (restricted to two-body problems), or correction via multiple-shooting schemes. These come at some non-trivial computational cost that must be paid each time a new trajectory is devised. The gravity of this issue can be partially attenuated by smart choice of coordinates (e.g. curvilinear, or orbit element differences), but it cannot be ignored.

\begin{figure}[h!]
\centering
\includegraphics[]{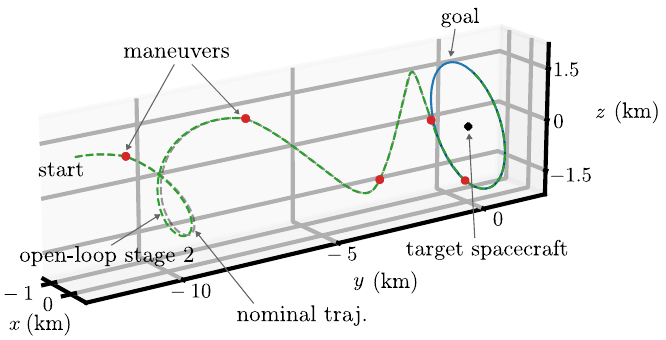}
\caption{Executed open-loop nonlinear guidance, nonlinear dynamics (rendezvous ex. 1)}
\label{fig:LEOGuidance_Plot3}
\end{figure}

In Fig.~\ref{fig:LEOGuidance_Plot3}, we show the results of Stage 2. This correction is achieved in $\sim$0.01s with 2 iterations of the Newton solver scheme ($\sim$0.005s per iteration) -- making use of analytic monomial expressions and the pre-saved $\Psi_{j}(t)$ evaluated at the critical times. The resulting nominal trajectory is again simulated in open-loop with the nonlinear dynamics, but now satisfies the final conditions to a high degree of accuracy (0.37\% error in final along-track position, $<0.1$\% error in all other states). In our tests, we find that this two-stage guidance scheme is quite stable, even when the linearized initial guess (Stage 1) is highly erroneous. 

\subsection{Nonlinear spacecraft rendezvous in LEO -- Sequential convexification-based optimal guidance}
\begin{figure}[h!]
\centering
\includegraphics[]{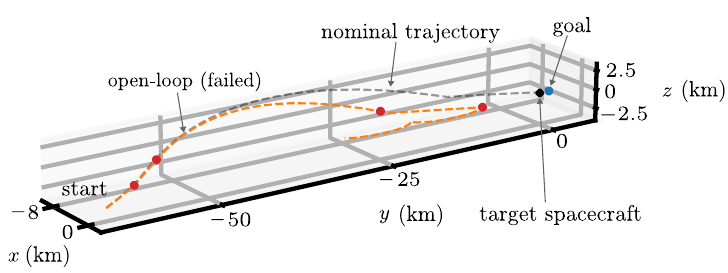}
\caption{Executed open-loop linear guidance, nonlinear dynamics (rendezvous ex. 2a)}
\label{fig:LEOGuidance_Plot6}
\end{figure}

Now we provide a simple demonstration of the sequential convexification scheme introduced in Section III.B. For this example we consider a new LEO rendezvous case whose details are provided in Table~\ref{table:LEO3}. This example combines medium range along-track separation with some out-of-plane motion which must be regulated to achieve a final stationary along-track state with 1.5 km offset from the target spacecraft. As before, we first compute a linearized guidance solution (Stage 1) and then a nonlinear guidance solution (now via the SCP scheme). The linearized guidance solution, and its failed open-loop execution, are depicted in Fig. \ref{fig:LEOGuidance_Plot6}. This solution is computed in $\sim$0.05 s.

\begin{table}[h!]
\centering
\caption{Parameters for Rendezvous Example 2a (Min. sum of delta-V squared, $\mathcal{P}=2$)}
\label{table:LEO3}
\begin{tabular}{ll} 
\hline\hline
\textbf{Category} & \textbf{Values} \\ 
\hline
\textbf{Target spacecraft} & LEO orbit, $a = 6378$ km, $e = 0$, $T \approx 1.4$ hrs \\ 
\textbf{Control interval} & $t_{0} = 0.1T$, $t_{f} = 1.1T$, discretized into 4 points (burn times fixed) \\ 
\textbf{Relative state at $t = 0$} & $\bm{x}(0) =[-3666.7, -62000, -4000, -1.239, 7.437, 2.479]$ (units: m, m/s) \\ 
\textbf{Relative state at $t_{f}$} & $\bm{x}(t_{f}) = [0, 1500, 0, 0, 0, 0]$ \\
\textbf{STT info} & Reused from Example 1 \\
\hline
\textbf{Linear guidance} & \\
Burn indices & $[0,12,64,99]$ \\ 
$\Delta v$ & 10.04 m/s \\ 
Runtime & 0.052s \\ 
Open-loop error & $>$10 km range, trajectory failed \\
\hline
\textbf{SCP guidance} & \\
$\Delta v$ & 10.82 m/s \\ 
Runtime & 0.1s \\ 
Open-loop $t_{f}$ error & 0.0529 km, 6.0 cm/s \\ 
\hline
\textbf{SCP parameters} & \\
Trust region & Fixed, $d = 3$ \\ 
Slack penalty weight & $w = 20$ \\ 
Method & ECOS \\ 
Convergence condition & $\|\tilde{\bm{X}}\| < 10^{-4}$ (converged in 5 iterations) \\
\hline\hline
\end{tabular}
\end{table}

\begin{figure}[h!]
\centering
\includegraphics[]{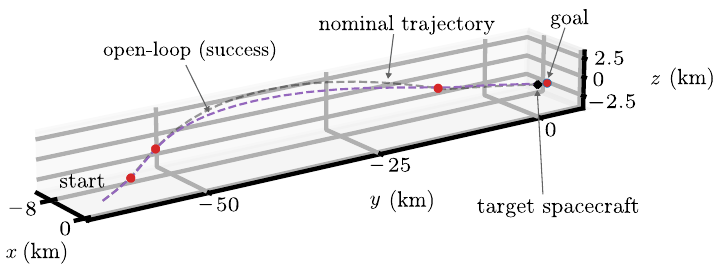}
\caption{Executed open-loop SCP-based guidance, nonlinear dynamics (rendezvous ex. 2a)}
\label{fig:LEOGuidance_Plot7}
\end{figure}

\begin{figure}[h!]
\centering
\includegraphics[]{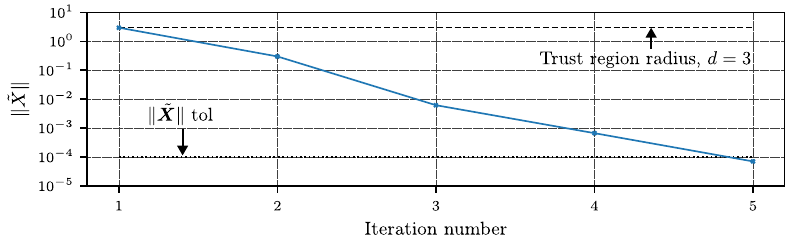}
\caption{SCP norm of non-slack free vars. vs. iter no. (rendezvous ex. 2a)}
\label{fig:SCP1_Xtilde}
\end{figure}

\begin{figure}[h!]
\centering
\includegraphics[]{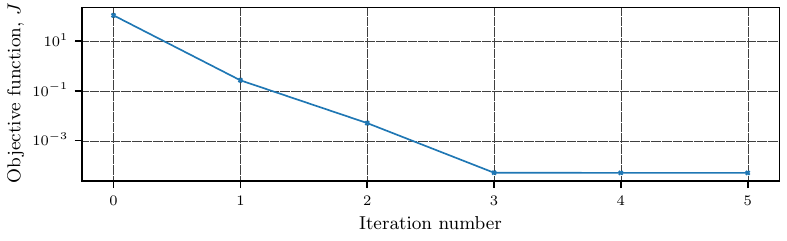}
\caption{SCP total cost vs. iter no. (rendezvous ex. 2a)}
\label{fig:SCP1_J}
\end{figure}

Figure \ref{fig:LEOGuidance_Plot7} gives the nonlinear guidance solution obtained via the SCP framework, minimizing sum of delta-V norm squared, thus $\mathcal{P}=2$, where the burn nodes and delta-V vectors are allowed to change via implementation of the scheme discussed in Section III.B. To facilitate the most rapid solution we inherit the optimal burn times of the linearized guidance solution, which minimizes the dimensionality of decision variables (although the SCP framework allows for the times to be optimized as well). The SCP solution, which is generated using the linearized guidance solution as an initial guess, converges in 5 iterations in a total runtime of $\sim$0.1 s, with a runtime per iteration of $\sim$0.02 s. The total solve time (linear SOCP + nonlinear SCP) to generate this nonlinear guidance solution was thus $\sim$ 0.15 s. 
The resulting guidance solution requires 10.82 m/s of total delta-V, compared to the linearized prediction of 10 m/s, and a two-stage guidance requirement of 11.2 m/s. In open-loop execution with the nonlinear dynamics, the guidance solution results in very modest final state error (as shown in Table 2), reflecting the limits of accuracy of the map $\Psi_{j}$ for order $j = 3$ with the monomials defined in the Cartesian coordinates.
For completeness, Figs.~\ref{fig:SCP1_Xtilde} and \ref{fig:SCP1_J} give the norm of the (non-slack) free variables, $\|\tilde{\bm{X}}\|$, and the total cost $J$ vs. iteration number. The steady drop in the norm of the free variables per iteration, as well as decreasing cost, correspond with our expectations of nominal behavior of the SCP scheme converging to a feasible minimum. For this problem, the norm of the vector of slack variables, $\|\bm{S}\|$, never exceeded $10^{-6}$.

\begin{figure}[h!]
\centering
\includegraphics[]{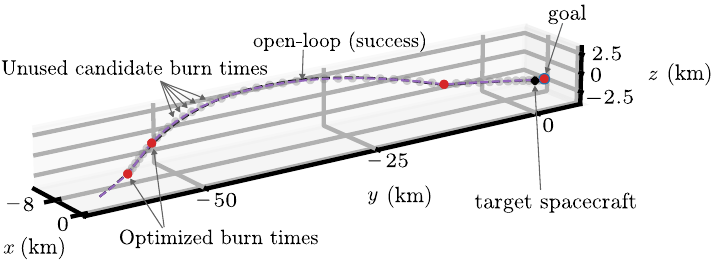}
\caption{Executed open-loop SCP-based guidance, nonlinear dynamics (rendezvous ex. 2b)}
\label{fig:LEOGuidance_Plot8}
\end{figure}

For this same problem, we also compute via SCP the $\mathcal{P}=1$ delta-V minimizing solution, and in addition we allow for the burn times to be selected from 100 candidate discretization points. The details for the new solver scheme are given in Table~\ref{table:LEO3b}. The final SCP solution still has only 4 maneuvers, but the timing of the inner two is shifted slightly. The other discretization times are automatically identified (as a by-product of the geometrically optimal path on $\mathcal{C}^{(N,j)}$) to be sub-optimal burn times, so the optimal result gives only 4 non-null jumps in the $\bm{c}_{j}(t_{i})$ across the 100 discretizations of the path. Note the increase in runtime to 2.5s as a result of the twenty-five fold expansion of the dimensionality of free variables for this example. The nominal delta-V is also slightly reduced. Figure~\ref{fig:LEOGuidance_Plot8} gives the new guidance solution, with the SCP solution in purple, the discretization of the problem shown by the gray dots (showing only 50 of the 100 points for clarity), and the selected burn nodes shown as red dots.
\begin{table}[h!]
\centering
\caption{Parameters for Rendezvous Example 2b (Min. sum of delta-V, $\mathcal{P}=1$)}
\label{table:LEO3b}
\begin{tabular}{ll} 
\hline\hline
\textbf{Category} & \textbf{Values} \\ 
\hline
\textbf{Target spacecraft} & LEO orbit, $a = 6378$ km, $e = 0$, $T \approx 1.4$ hrs \\ 
\textbf{Control interval} & $t_{0} = 0.1T$, $t_{f} = 1.1T$, discretized into 100 points (burn times free) \\ 
\hline
\textbf{Linear vs. SCP} & \\
Burn indices (Linear) & $[0,12,64,99]$ \\ 
Burn indices (SCP solution) & $[0,13,65,99]$ \\ 
\hline
\textbf{SCP guidance} & \\
$\Delta v$ & 10.73 m/s \\ 
Runtime & 2.5s \\ 
Open-loop $t_{f}$ error & 0.108 km, 14.6 cm/s \\
\hline
\textbf{SCP parameters} & \\
Trust region & Fixed, $d = 10$ \\ 
Slack penalty weight & $w = 5$ \\ 
Method & ECOS \\ 
Convergence condition & $\|\tilde{\bm{X}}\| < 10^{-4}$, (converged in 6 iterations) \\
\hline\hline
\end{tabular}
\end{table}

All guidance examples shown so far were generated using the same pre-computed third order nonlinear expansion, and represent a small sample of the trade space that can be explored leveraging that same data. The two-stage nonlinear guidance and SCP-based nonlinear guidance implementations all perform without the need for any real-time numerical integration, and their robustness to starting with a poor linearization-based initial guess is encouraging and interesting. Both offer nonlinear correction in a comparable amount of time that was required to generate the linearized initial guess, but of course the SCP-based implementation comes with the added delta-V minimizing benefit. Note that our implementation in Python is not optimized for speed but for proof of concept, and a CVXPYgen implementation of the SCP could potentially improve performance by removing significant overhead from successive CVXPY calls \citep{cvxpygen2022}. 
\subsection{Long range rendezvous guidance via spherical coordinate monomials}
In solving the trajectory optimization problem, any desired working coordinates can be used in the construction of the $\bm{c}_{j}$ and the computation of their associated fundamental solution matrix $\Psi_{j}$. For orbital mechanics problems, it is well-known that curvilinear/spherical coordinates yield a lower level of dynamical nonlinearity than the Cartesian coordinates \citep{Willis2ndOrder, junkAdven}. In other words, for a given order $j$, the $\Psi_{j}$ computed in spherical coordinates will have a greater physical region of validity than its counterpart in Cartesian coordinates. In this section we modify the previous guidance implementation, instead now making use of a spherical coordinate overparameterized monomial representation. The result is a solver with a much greater region of validity ($\mathcal{O}(1000)$ km in LEO for order $j=3$, as opposed to $\mathcal{O}(100)$ km for the third-order Cartesian representation).

The relative position can be developed in relative spherical coordinates $\delta r, \delta\theta, \delta\phi$, which relate to the Cartesian relative position coordinates as:
\begin{equation}
\label{Coords1}
\begin{split}
x & \ = \left(r + \delta r\right)\cos{\delta\theta}\cos{\delta\psi} - r \\
y & \ = \left(r + \delta r\right)\sin{\delta\theta}\cos{\delta\phi} \\
z & \ = \left(r + \delta r\right)\sin{\delta\phi}
\end{split}
\end{equation}
where $r$ is the instantaneous target orbit radius, $r_{c} = r + \delta r$ is the instantaneous chaser orbit radius, $\delta\theta$ is an Earth-centered angular offset measured in the plane of the target orbit, and $\delta\phi$ is the out of plane angular offset. The derivatives of these quantities,  $\delta\dot{r}, \delta\dot{\theta}, \delta\dot{\phi}$, constitute the remainder of the state description, and the Cartesian velocities relate to the spherical relative states as below: 
\begin{equation}
\label{Coords2}
\begin{split}
\dot{x} & \ = (\dot{r} + \delta\dot{r})\cos{\delta\phi}\cos{\delta\theta} - (r + \delta r)(\delta\dot{\phi}\sin{\delta\phi}\cos{\delta\theta} + \delta\dot{\theta}\cos{\delta\phi}\sin{\delta\theta}) - \dot{r} \\
\dot{y} & \ = (\dot{r} + \delta\dot{r})\cos{\delta\phi}\sin{\delta\theta} - (r + \delta r)(\delta\dot{\phi}\sin{\delta\phi}\sin{\delta\theta} - \delta\dot{\theta}\cos{\delta\phi}\cos{\delta\theta}) \\
\dot{z} & \ = (\dot{r} + \delta\dot{r})\sin{\delta\phi} + (r + \delta r)\delta\dot{\phi}\cos{\delta\phi}
\end{split}
\end{equation}
Using the notation of Section II.E.2, these define new working coordinates, $\bm{\eta} = (\delta r, \delta\theta, \delta\phi, \delta\dot{r}, \delta\dot{\theta}, \delta\dot{\phi})^{\top}$. See e.g. \cite{willis2019secondorder} for the full inverse mapping between Cartesian and spherical coordinates. For convenience we use the normalized spherical relative coordinates $\overline{\bm{\eta}} = \left(\delta\rho, \delta\theta, \delta\phi, \delta\rho', \delta\theta', \delta\phi'\right)^{\top}$. The radial coordinates are normalized by $\delta\rho = \delta r/a$, and we normalize time as $\zeta = nt$, thus $( \ )^{'} = \frac{\text{d}}{\text{d}\zeta}( \ ) = \frac{1}{n}\frac{\text{d}}{\text{d}t}( \ )$. 

We develop code to compute new maps $\Psi_{j}(t)$ that contain the nonlinear fundamental solutions for the following 3D (6 state) normalized target-centered spherical coordinate nonlinear dynamics of Keplerian relative motion:
\begin{subequations}
\label{SphNondim1}
\begin{align}
\delta\rho'' & \ = \frac{1}{\rho^{2}} - \frac{1}{\left(\rho + \delta\rho\right)^{2}} - \rho\theta'^{2} + (\rho + \delta\rho)\left(\delta\phi'^{2} + (\theta' + \delta\theta')^{2}\cos^{2}{\delta\psi}\right) \\
\delta\theta'' & \ = 2\frac{\rho'\theta'}{\rho} - 2\frac{\rho' + \delta\rho'}{\rho + \delta\rho}\left(\theta' + \delta\theta'\right) + 2\left(\theta' + \delta\theta'\right)\delta\phi'\tan{\delta\phi} \\
\delta\phi'' & \ = -2\frac{\rho' + \delta\rho'}{\rho + \delta\rho}\delta\phi' - \left(\theta' + \delta\theta'\right)^{2}\cos{\delta\phi}\sin{\delta\phi}
\end{align}
\end{subequations}
These are integrated in parallel with the normalized 2D (4 state) equations of the target orbit:
\begin{subequations}
\label{SphNondim2}
\begin{align}
\rho'' = & \ -\frac{1}{\rho^{2}} + \rho\theta'^{2} \\
\theta'' = & \ -2\frac{\rho'\theta'}{\rho}
\end{align}
\end{subequations}

where $\rho$ is the target orbital radius and $\theta$ its argument of latitude. This normalization scheme allows the same $\Psi_{j}$ to be reused for rendezvous guidance for all target orbits of a given eccentricity, regardless of the target semimajor axis. This enables all circular Keplerian rendezvous guidance problems of a given maximum normalized duration, discretization, and scale to be handled using the same pre-computed information. Using Eq.~\eqref{Coords2} and the normalization definition, and denoting the transformation from normalized spherical states to dimensional Cartesian velocity as $\bm{v} = \bm{g}_{v}(\overline{\bm{\eta}})$, the delta-V resulting from a given impulsive change in the normalized spherical coordinates can be defined:
\begin{equation}
\label{DvSph}
\Delta\bm{v}(t_{i}) = \bm{g}_{v}(\overline{\eta}(t_{i}^{+})) - \bm{g}_{v}(\overline{\eta}(t_{i}^{-}))
\end{equation}
with $t_{i}^{-}$ and $t_{i}^{+}$ denoting the state just before and after the maneuver. This is a complex \textit{nonlinear} function of the normalized spherical coordinates, but it admits (approximately) the following linear relationship to the overparameterized spherical monomial states, computed using differential algebra:
\begin{equation}
\label{DvSph2}
\Delta\bm{v}(t_{i}) = \Gamma_{v,j}(t_{i})\left( \bm{c}_{j}^{(\overline{\eta})}(t_{i}) - \bm{c}_{j}^{(\overline{\eta})}(t_{i-1}) \right)
\end{equation}
This transformation renders delta-V norm and norm-squared cost functions convex in $\bm{c}_{j}^{(\overline{\eta})}$. As for problem constraints, we require that impulsive maneuvers change only the velocity components $\delta\rho', \delta\theta', \delta\phi'$ without affecting the position components $\delta\rho, \delta\theta, \delta\phi$. As a result, the kinematic constraints are otherwise still handled by the $\Psi_{j}$ (but now developed for the spherical monomials), and a spherical coordinate equivalent of the SCP scheme defined in Eq.~\eqref{TrajOpt6} modifies only the cost terms, with a substitution of $\Gamma_{v,j}(t_{i})$ for the $\Psi_{v,j}(t_{i})$.

Table~\ref{table:Sph1} gives some info about the maps $\Psi_{j}$ computed for the spherical monomials using differential algebra methods. We compute up to order $j=4$, and note at order $j=2$ that we recover numerically the analytic finding of \cite{willis2019secondorder} that of 21 possible second-order fundamental solutions (reflecting the 21 different quadratic monomials), only 15 are nonzero. This is opposed to the Cartesian coordinate solution, for which 19 are nonzero. Thus the spherical coordinate parameterization is more compact. This property extends to order $j=4$, at which point of the 209 possible fundamental solutions, corresponding with all unique monomials from orders 1 to 4, 83 are zero. We also determine that the computation speed of the map is not greatly affected by eccentricity, but more eccentric cases require a finer discretization in the integration of Eqs.~\eqref{SphNondim1} and Eq.~\eqref{SphNondim2} to keep the same level of accuracy. Due to the normalization, a computed map of given eccentricity applies for any target semimajor axis and orbit period, about any planet. Equivalent maps could be constructed for different eccentricities, or for altogether different dynamical regimes, but we continue with the classic circular orbit case $e=0$ studied most frequently in the spacecraft rendezvous community. We write a new SCP solver using the normalized spherical coordinates as our working coordinates for the monomial basis. The normalization also standardizes the choice of (normalized) trust region $d$ vs. order $j$, so that spacecraft rendezvous expertise is no longer required in choosing this value. 
\begin{table}[h!]
\centering
\caption{Developing the Fundamental Solution Matrix for Various Cases}
\label{table:Sph1}
\begin{tabular}{ll} 
\hline\hline
\textbf{Category} & \textbf{Values} \\ 
\hline
\textbf{Discretization info} & $\Delta \zeta = 4\pi$ (2 target orbits), discretized into 400 points \\ 
\hline
\textbf{Circular orbit ($e = 0$)} & \\
Compute time (order $j$) & 3.2s (4), 1.4s (3), 0.6s (2), 0.1s (1) \\ 
Data size (order $j=4$) & 2.96 MB \\ 
Number of empty columns of $\Psi_{j}$ for $j=4$ & 83 \\ 
\hline
\textbf{Eccentric orbit ($e = 0.2$)} & \\
Compute time (order $j$) & 3.2s (4), 1.4s (3), 0.6s (2), 0.3s (1) \\ 
Data size (order $j=4$) & 3.04 MB \\ 
Number of empty columns of $\Psi_{j}$ for $j=4$ & 83 \\ 
\hline\hline
\end{tabular}
\end{table}

\begin{figure}[h!]
\centering
\includegraphics[width=6in]{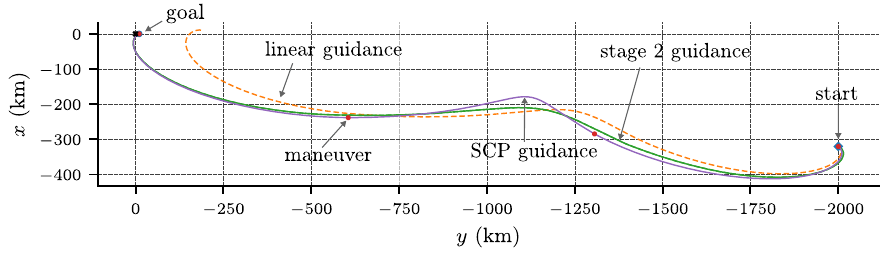}
\caption{Executed open-loop guidance, SCP with inherited burn times (rendezvous ex. 3a)}
\label{fig:LEOGuidance_LR1}
\end{figure}

\begin{figure}[h!]
\centering
\includegraphics[width=6in]{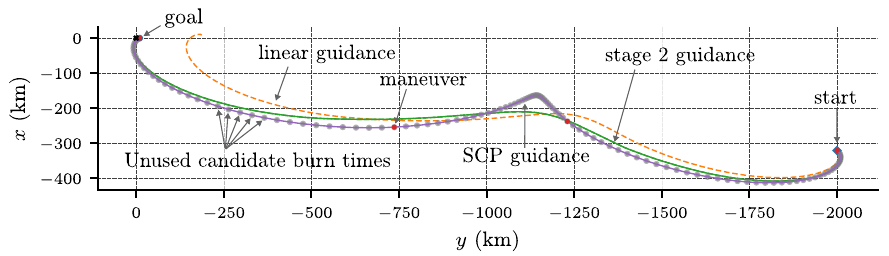}
\caption{Executed open-loop guidance, SCP with free burn times (rendezvous ex. 3b)}
\label{fig:LEOGuidance_LR2}
\end{figure}

Consider the long-range case described in Table~\ref{table:ULR1}. The chaser spacecraft is initialized in an orbit of a similar period but phased over 2000 km behind the target spacecraft, with some additional out-of-plane motion that must be regulated. The control problem is to reduce the inter-spacecraft separation to a purely along-track separation of 10 km, whereupon terminal rendezvous procedures can later be initiated. Of the 400 saved times for the map $\Psi_{j}$, this data is down-sampled by a factor of three to define a control problem with a reasonable number of candidate burn times. To cover the large distance in just under 2.5 hours, the optimal solution will use much more delta-V than the earlier examples. We provide two solutions. The first solution, depicted in Fig.~\ref{fig:LEOGuidance_LR1}, uses the two-stage guidance scheme to produce a dynamically feasible initial guess of the optimal transfer, which is then corrected using the SCP scheme. The linearly predicted optimal burn times are inherited, but the location and delta-V vector of these burn nodes is adjusted via the SCP. For the second solution, given by Fig.~\ref{fig:LEOGuidance_LR2}, the two-stage guidance solution is corrected with the burn times left free, and the SCP scheme changes the timing of the burns slightly (each $<$5 discretized times away from all nominal burns). The result is a slightly different shape to the trajectory: observe the ``kink" during which the relative motion trajectory greatly slows down, as evidenced by the clustering of discretization points which are chosen uniformly in time for this example. 
\begin{table}[h!]
\centering
\caption{Simulation Parameters for Rendezvous Example 3 (Min. sum of delta-V, $\mathcal{P}=1$)}
\label{table:ULR1}
\begin{tabular}{ll} \hline\hline
Parameter                         & Values \\
\cline{1-1}  \cline{2-2}  \vspace{-0.1cm}
Target spacecraft           & LEO orbit, $a = 6378$ km, $e = 0$, $T \approx 1.4$ hrs                                            \\
Control interval                  & $t_{0} = 0.05T$, $t_{f} = 1.80T$, discretized into 118 points                       \\
Relative state, $t = 0$ & $\bm{x}(0) =[-320.4, -2000, 70, 0, -5n, 10n]$ (units: km, km/s), $n=2\pi/T$                                                      \\
Relative state, $t_{f}$  & $\bm{x}(t_{f}) = [0, -10, 0, 0, 0, 0]$ \\
 \hline\hline
\end{tabular}
\end{table}

\begin{table}[h!]
\centering
\caption{Guidance Results for Rendezvous Example 3}
\label{table:ULR2}
\begin{tabular}{ll} 
\hline\hline
\textbf{Category} & \textbf{Values} \\ 
\hline
\textbf{Two-stage guidance} & \\
Burn indices & $[0, 30, 40, 83, 116]$ \\ 
$\Delta v$ & 243.00 m/s \\ 
Total runtime & 0.128s \\ 
Open-loop error (pos/vel norm) & 0.645 km, 0.356 m/s \\ 
No. iterations needed & 3 (which took 0.054s) \\ 
\hline
\textbf{SCP parameters} & \\
Trust region & Fixed, $d = 0.1$ \\ 
Slack penalty weight & $w = 10^{3}$ \\ 
Method & ECOS \\ 
Convergence condition & $\|\tilde{\bm{X}}\| < 10^{-7}$ (note: problem normalized) \\ 
\hline
\textbf{SCP (fixed burn times)} & \\
Burn indices & $[0, 40, 83, 116]$ \\ 
$\Delta v$ & 226.31 m/s \\ 
Runtime & 0.296s \\ 
Open-loop end error (pos/vel norm), $j=4$ & 0.038 km, 0.193 m/s \\ 
No. iterations needed & 6 \\ 
\hline
\textbf{SCP (free burn times)} & \\
Burn indices & $[0, 44, 81, 117]$ \\ 
$\Delta v$ & 222.22 m/s \\ 
Runtime & 3.89s \\ 
Open-loop end error (pos/vel norm), $j=4$ & 0.0857 km, 0.332 m/s \\ 
No. iterations needed & 10 \\ 
\hline\hline
\end{tabular}
\end{table}

Table~\ref{table:ULR2} gives the guidance solution info, as well as runtime and open-loop execution error for the $j=4$ solutions. In comparison to these results, note that $j = 3$ produced a highly similar guidance solution, but also resulted in final state targeting range error of $\sim$1 km. While this is an error on the scale of $\sim$0.05\% in terms of the initial range, it is a $\sim$10\% error in the much smaller goal range, necessitating the extra order to bring this error below 1\%. The $j=4$ solution required only $\sim$47\% more compute time than $j=3$. Note in both SCP solutions, the nonlinear optimal solution removes the need for the linearly predicted optimal burn at index 30 entirely. However the two-stage guidance already provides a decent solution, with SCP results only reducing total delta-V on the order of 10\%. In comparing the SCP solution which inherits the linearly optimal burn times to the solution with burn times free, note that the burn times are shifted slightly, the delta-V is decreased by another $\sim$2\%, but the runtime is over 10 times longer due to a similar expansion in the number of decision variables. Overall, this demonstration illustrates a few concepts: 1) the importance of choice of working coordinates for the monomials and the superiority of spherical coordinates over Cartesian, 2) the usefulness of problem normalization, and 3) the increase in runtime with order $j$, which is much less than the growth of the dimensionality of the monomials, $K_{j}$.

\subsection{Nonlinear state constraints}
Finally, we apply a simple set of constraints to modify the prior example and show the property of embedding nonlinear constraints. We consider a series of constraints on the inter-spacecraft range as follows:
\begin{equation}
\label{rangeConstr1}
\varrho_{\text{min}} = \left \{ \begin{array}{ll} 1915\text{ km} & t\leq 0.4T  \\ 1405\text{ km} & 0.4T < t\leq0.75T \\ 318\text{ km} & 0.75T < t\leq1.4T \\ 64\text{ km} & 1.4T < t\leq1.6T \\ 7.5\text{ km} & t>1.6T \end{array}\right.
\end{equation}
These timed constraints emulate a simple active safety requirement, wherein the chaser spacecraft should not approach the target in a manner more aggressive than the specified approach stages. The inter-spacecraft range $\varrho$ is related to the instantaneous spherical relative position coordinates $\delta r, \delta\theta, \delta\phi$ and the target orbital radius $r$ by  the following nonlinear expression:
\begin{equation}
\label{rangeConstr2}
\begin{split}
\varrho & \ = \sqrt{\delta r^{2} + 2r\delta r + 2r^{2} - 2r\left(r + \delta r\right)\cos{\delta\phi}\cos{\delta\theta}} \\
& \ = a\sqrt{\delta \rho^{2} + 2\rho\delta \rho + 2\rho^{2} - 2\rho\left(\rho + \delta \rho\right)\cos{\delta\phi}\cos{\delta\theta}}
\end{split}
\end{equation}

Equivalently we can define the constraint $\varrho^{2}(t) \geq \varrho_{\text{min}}^{2}(t)$ to avoid the needless square root for an always-positive function. Furthermore we can normalize the constraint as $\overline{\varrho}^{2}(t) \geq \overline{\varrho}_{\text{min}}^{2}(t)$, with $\varrho = a\overline{\varrho}$. This constraint is nonlinear in the spherical relative coordinates, but via differential algebra we can impose this nonlinear constraint function as a linear constraint function on the $\bm{c}_{j}^{(\overline{\eta})}(t_{i}) \ \forall t_{i}$:
\begin{subequations}
\label{rangeConstr3}
\begin{align}
h(\overline{\bm{\eta}}(t_{i})) & \ = \overline{\varrho}^{2}(t_{i}) \\
h(\overline{\bm{\eta}}(t_{i})) & \ \geq \overline{\varrho}_{\text{min}}^{2}(t_{i}) \\
h(\overline{\bm{\eta}}(t_{i}), t_{k}) & \ = \bm{\gamma}_{h,j}^{\top}(t_{k})\bm{c}_{j}^{(\overline{\eta})}(t_{i}) \\
\bm{\gamma}_{h,j}^{\top}(t_{k})\bm{c}_{j}^{(\overline{\eta})}(t_{i}) & \ \geq \overline{\varrho}_{\text{min}}^{2}(t_{i})
\end{align}
\end{subequations}
where $\bm{\gamma}_{h,j}(t_{k})$ is computed from the normalized version of Eq.~\eqref{rangeConstr2} using differential algebra, similarly to $\Psi_{j}$ and $\Gamma_{v,j}$. This function provides the mapping from the monomial expansion $\bm{E}_{j}(\overline{\bm{\eta}}(0))$ to range at time $t_{k}$ via the $j$\textsuperscript{th}-order truncation of the flow of the natural dynamics. Thus it is also possible via this formulation to enforce passive safety constraints, wherein states at time $t_{i}$ are constrained based on their free-drift implications at times $t_{k} \geq t_{i}$, although this was not explored, and we enforce only instantaneous range at all trajectory nodes i.e. $t_{k} = t_{i}$.

Similarly to the kinematic constraints, Eq.~\eqref{rangeConstr3} is linear in the $\bm{c}_{j}^{(\overline{\eta})}$, thus the same sequential convexification procedure extends naturally to accommodate this extra constraint, given below as an addendum to Eq.~\eqref{TrajOpt6}:
\begin{equation}
\label{rangeConstr4}
\bm{\gamma}_{h,j}^{\top}(t_{i})\left(\bm{c}_{j}^{(\overline{\eta})^{\dagger}}(t_{i}) +  \left.\frac{\partial\bm{c}_{j}^{(\overline{\eta})}}{\partial\bm{c}_{1}^{(\overline{\eta})}}\right\vert_{i}\delta\bm{c}_{1}^{(\overline{\eta})}(t_{i}) \right) + s_{\text{ineq},i} \geq \overline{\varrho}_{\text{min}}^{2}(t_{i}), \ \ i = 1\text{:}K
\end{equation}
where the prime notation is reused to denote the  result of a prior iteration, and this constraint results in the introduction of an extra $K$ slack variables, expanding the dimensionality of the free variable vector $\bm{X}$ in the convex solver. These extra slack variables must also be penalized, necessitating the addition of the following cost term:
\begin{equation}
\label{rangeCost}
\delta J_{\text{ineq}} = \bm{s}_{\text{ineq}}^{\top}P_{w,\text{ineq}}\bm{s}_{\text{ineq}}
\end{equation}
where $\bm{s}_{\text{ineq}}$ is a vector of just the scalar inequality slack variables $s_{\text{ineq},i}$, and $P_{w,\text{ineq}}$ is a diagonal positive-definite matrix whose diagonal elements can be chosen so that all slack variables are given equal treatment.

\begin{figure}[h!]
\centering
\includegraphics[width=6in]{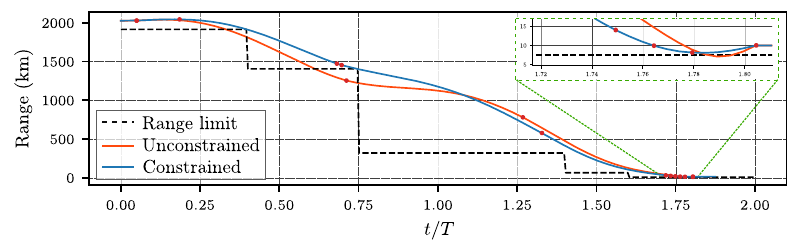}
\caption{Inter-spacecraft range, unconstrained and constrained (rendezvous ex. 3c)}
\label{fig:LEOGuidance_constrained}
\end{figure}

Applying the modified SCP scheme with path constraints as specified in Eq.~\eqref{rangeConstr1}, the inter-spacecraft range is given in Fig.~\ref{fig:LEOGuidance_constrained}. This figure depicts the time-dependent range constraints, showing the unconstrained result equivalent to Example 3(b), and the new result which satisfies the range constraints. We note the interesting addition of multiple burns at the end of the trajectory to avoid violation of the final range constraint of $\varrho_{\text{min}} = 7.5$ km, which was previously violated. The satisfactory constrained trajectory is obtained with a total delta-V of $239.66$ m/s, a 7.8\% increase on the unconstrained case. The new scheme converges in 9 iterations with a solve time of 5.86s, achieving final position and velocity error norms of 0.092 km and 0.155 m/s, of similar targeting performance to the prior unconstrained $j=4$ results. Note that additional modifications could be made to the solver to avoid returning sequences of small maneuvers -- this is a consequence of our problem formulation and our simple delta-V minimizing cost function. The adjacent burns can usually be combined for fairly small delta-V penalty. Overall this example illustrates the inclusion of nonlinear state constraints into our new methodology, and we reserve more refined solver development to future work.

\section{Conclusions}
This paper introduces a method of rapid guidance for nonlinear systems leveraging overparameterized monomial coordinates and fundamental solution expansions. In essence, we trade nonlinearities and non-convexities of the original problem for a single non-convexity of our own design, posing the problem on the non-Euclidean surface $\mathcal{C}^{(N,j)}$ of known structure.  By this methodology, trajectory optimization problems can be solved as path planning problems with very low real-time computational footprint. The solutions inherit the accuracy of the nonlinear fundamental solution expansions encoded by the the special nonlinear fundamental solution matrix $\Psi_{j}(t)$. This is in contrast to collocation-based SCP schemes which have no guarantee of physical feasibility. We devise a simple two-stage (linear predict, nonlinear correct) approach and a simple sequential convexification scheme to demonstrate the usefulness of the methodology. The spacecraft relative motion and rendezvous problem (and relatedly, orbit regulation) are particularly well-suited for this approach because the reference trajectory needed for the nonlinear solution expansions is simply the target spacecraft orbit. We present several examples where our methodology is used to rapidly compute guidance solutions that are effective when propagated in a numerical model of nonlinear problem dynamics. 

In general, like linear methods, this method is still only ``locally" valid, though admitting a much larger domain of validity than linearization. For this reason, application to more general orbit transfer problems requires the computation of reasonable reference trajectories for the transfer, which are often not obvious, particularly in a three- or four-body context. Also, in this work, we consider only reference trajectories that are uncontrolled natural solutions to the nonlinear dynamics. This is not particularly limiting for the chosen example of rendezvous and station-keeping, where the choice of reference is clear.  Furthermore, only an impulsive maneuver guidance approach is developed, to entirely avoid the use of collocation and interpolation schemes. It is also worth noting that, depending on order of nonlinearity $j$, the time to compute the fundamental solutions numerically can be non-negligible. We envision a GNC implementation where these solutions are computed hours, days, or more in advance of when needed in a mission. In the context of our chosen example problem of spacecraft rendezvous, incorporation of analytic or numerical fundamental solution expansions keeps the total computational footprint extremely low (solve times of $\mathcal{O}(1\text{s})$ or shorter for our examples). Overall, this work establishes a foothold in a new research direction with high promise for fast and reliable guidance suitable for embedded systems in spaceflight. We expect that the themes and techniques explored will find further use by many other researchers interested in computationally lean optimization.

\section*{Disclaimer} 
Funded by the European Union. Views and opinions expressed are however those of the authors only and do not necessarily reflect those of the European Union or the European Research Executive Agency. Neither the European Union nor the granting authority can be held responsible for them. 

\section*{Code and Data Availability}
Data and links to code to reproduce select examples from this paper will be made available via the open-access digital repository Zenodo at \url{https://zenodo.org/communities/faast-msca/about}.

\section*{Acknowledgements} 
This work was performed as part of the MSCA project ``Facilitating Autonomy in Astrodynamics for Spacecraft Technology", Grant agreement ID: 101063274. We gratefully acknowledge technical and administrative support from the MSCA hosting organization: Department of Aerospace Science and Technology (DAER), Politecnico di Milano, Italy.

\bibliography{references.bib}   

\appendix

\section{Computerized generation and manipulation of monomial representations}
To implement our methodology, one needs a well-ordered and unambiguous way of constructing and manipulating monomial sequences. We use a positional power representation of monomials and we represent sequences of monomials in 2D arrays. The easiest way to explain this is with a full end-to-end example, but \cite{AlgebraicManip, JorbaNormalForms} give more context on similar techniques. Here we drop the initial condition notation ``$x_{i}(0)^{j}$" in favor of the more compact $x_{i}^{j}$. 

Consider a 3-state system with states $x_{1}$, $x_{2}$, $x_{3}$. The array representation of the monomial $x_{1}^{a}x_{2}^{b}x_{3}^{c}$ for integer powers $\{a,b,c\} \geq 0$ is $[a,b,c]$. These are stacked row-wise to reproduce sequences of monomials as needed. For example, we represent $\bm{c}_{1} = \bm{x} = (x_{1},x_{2},x_{3})^{\top}$ as the following array:

\begin{equation}
\label{arrC1}
\text{arr}\left(\bm{c}_{1}\right) = \begin{bmatrix} 1 & 0 & 0  \\ 0 & 1 & 0 \\ 0 & 0 & 1   \end{bmatrix}
\end{equation}

The following pseudocode constructs all rows of a desired $\text{arr}\left(\bm{c}_{j}\right)$ for a system with $N$ states by initializing $\text{arr}\left(\bm{c}_{1}\right)$ and then, order-by-order, appending the array representation (itself row-by-row) for higher-order terms up to and including nonlinearity order $j$. We take the array representation of a prior order (``base\_block") and generate all unique monomial power permutations for a new order (``new\_block"), then step to the next order with base\_block $\leftarrow$ new\_block.

\vspace{12pt}
\begin{breakablealgorithm}
  \begin{algorithmic}[1]
    \Procedure{Build\_cj\_Arr}{$N, \ j$}\Comment{Construct, $\text{arr}(\bm{c}_{j})$, the $K_{j}\times N$ array representation of $\bm{c}_{j}$}
    \State $K_{j} \gets \text{compute\_Kj}(N,j)$\Comment{Apply Eq.~\eqref{WhatIsKj1} to get total number of monomials for orders 1 through $j$}
    \State $\text{arr}(\bm{c}_{j}) \gets I_{N\times N}$\Comment{Initialize $\text{arr}(\bm{c}_{j})$ at order $j=1$ as the trivial $N\times N$ representation for $\bm{c}_{1}$}
    \State mult\_mat $ \gets I_{N\times N}$\Comment{Multiplication by state $x_{q}$ will be achieved via addition of $q$\textsuperscript{th} row of mult\_mat}
   \State  base\_block $\gets \text{arr}(\bm{c}_{j})$\Comment{Initialize base\_block}
   \State new\_block $\gets$ [ ]\Comment{new\_block is initialized empty}
    \State 
    \For{$i_{1}$ from 2 through $j$}\Comment{This order}
        \State n\_rows $\gets$ base\_block.num\_rows\Comment{Number of rows in base\_block}
                \For{$i_{2}$ from 1 through $N$}\Comment{This state}
        \State term\_mult $\gets$ mult\_mat[$i_{2}$, :]\Comment{Multiplier for new monomial candidate}
              \For{$i_{3}$ from 1 through n\_rows}\Comment{This row of base\_block}
              	\State new\_row $\gets$ base\_block[$i_{3}$, :] + term\_mult\Comment{Make representation of new candidate monomial}
		\If{new\_row $\notin$ new\_block}
		\State new\_block.append(new\_row)\Comment{Append candidate monomial if it is new}
		\EndIf
              \EndFor
        \EndFor
    \State $\text{arr}(\bm{c}_{j})$.append(new\_block)\Comment{Append base\_block; new order done}
    \State base\_block $\gets$ new\_block\Comment{step}
    \EndFor
    \State \textbf{return} $\text{arr}(\bm{c}_{j})$\Comment{Return the completed array representation of $\bm{c}_{j}$}
    \EndProcedure
  \end{algorithmic}
\end{breakablealgorithm}
\vspace{12pt}

As an example output, we provide explicitly the third-order array representation of $\bm{c}_{j}$ for the case of $N = j = 3$ below, with sub-blocks of orders 1, 2, and 3 partitioned for clarity:

\begin{equation}
\label{arrC2}
\text{arr}\left(\bm{c}_{3}\right) = \begin{bmatrix} 1 & 0 & 0  \\ 0 & 1 & 0 \\ 0 & 0 & 1 
 \\ \hdashline[2pt/2pt] \\[-2.5ex] 2 & 0 & 0 \\ 1 & 1 & 0 \\ 1 & 0 & 1 \\ 0 & 2 & 0 \\ 0 & 1 & 1 \\ 0 & 0 & 2 
  \\ \hdashline[2pt/2pt] \\[-2.5ex] 3 & 0 & 0 \\ 2 & 1 & 0 \\ 2 & 0 & 1 \\ 1 & 2 & 0 \\ 1 & 1 & 1 \\ 1 & 0 & 2 \\ 0 & 3 & 0 \\ 0 & 2 & 1 \\ 0 & 1 & 2 \\ 0 & 0 & 3
 \end{bmatrix}
\end{equation}

\noindent where e.g. the fourth row (first second-order component) is $x_{1}^{2}$, the eighth row is $x_{2}x_{3}$, and the eleventh is $x_{1}^{2}x_{2}$. Note that for nonlinearities of order $q$, it can be shown that the resultant sub-block will have the following number of terms:

\begin{equation}
\label{arrC3}
\text{dim}\left ( \text{sub-block } q\right) = \begin{pmatrix} N + q - 1 \\ q \end{pmatrix}
\end{equation}

Adding up the number of rows of all sub-blocks to the max order $j$, we recover Eq.~\eqref{WhatIsKj1}. 

With a computerized array representation of $\bm{c}_{j}$, the column-wise arrangement of all fundamental solutions to make $\Psi_{j}$ is also specified. We defer the reader to our code for any more details in that regard. We now briefly discuss how the Jacobians $C_{j} = \frac{\partial\bm{c}_{j}}{\partial\bm{c}_{1}}$ are computed in a computerized fashion leveraging the array representation of $\bm{c}_{j}$. This Jacobian should be of size $K_{j} \times N$ and will be composed completely of monomials and scalars. To compute the elements of this it is easy to differentiate a row of $\text{arr}(\bm{c}_{j})$ with respect to a state by removing a one from the entry corresponding to the state being differentiated, and retaining the prior power as a new multiplying coefficient. For example, $\frac{\partial}{\partial x_{3}}\left(1\times [1, 0, 2]\right) = 2\times [1, 0, 1]$, and $\frac{\partial}{\partial x_{3}}\left(1\times [1, 1, 0]\right) = 0\times [\ldots]$. In this manner derivatives on the monomials are reduced to simple array operations, which renders a very lean Jacobian computation. Once an array representation is developed, it can be reused, and furthermore the Jacobian is a linear function of $\bm{c}_{j}$, which is computationally convenient.

\section{Computing the nonlinear fundamental solution expansion from STTs}
This is the first of three admissible methods for constructing $\Psi_{j}(t)$. \cite{ParkScheeresSTTs} introduces the state transition tensors (STTs) as a generalization of the more familiar state transition matrix. They map initial deviations of a system to deviations at some time $t$, and are characterized by the scalar expansions below, borrowing their notation:
\begin{subequations}
\label{STT1}
\begin{align}
\delta x_{i}(t) = \sum_{p=1}^{m}\frac{1}{p!}\Phi_{i,k_{1}\ldots k_{p}}\delta x_{k_{1}}^{0}\ldots\delta x_{k_{p}}^{0}, \ i = 1:N \\ \delta\dot{x}_{i}(t) = \sum_{p=1}^{m}\frac{1}{p!}f^{*}_{i,k_{1}\ldots k_{p}}\delta x_{k_{1}}^{0}\ldots\delta x_{k_{p}}^{0} \\ \delta\dot{x}_{i}(t) = \sum_{p=1}^{m}\frac{1}{p!}\dot{\Phi}_{i,k_{1}\ldots k_{p}}\delta x_{k_{1}}^{0}\ldots\delta x_{k_{p}}^{0}
\end{align}
\end{subequations}
where summation convention is used, $k_{j} \in \{1, \ldots, N\}$, and subscripts $k_{j}$ denote the $k_{j}$\textsuperscript{th} component of the state vector:
\begin{equation}
\label{STT2}
\Phi_{i,k_{1}\ldots k_{p}} = \frac{\partial^{p}x_{i}}{\partial x_{k_{1}}^{0}\ldots\partial x_{k_{p}}^{0}}
\end{equation}
\begin{equation}
\label{STT3}
f_{i,k_{1}\ldots k_{p}}^{*} = \left.\frac{\partial^{p}f_{i}}{\partial x_{k_{1}}^{0}\ldots\partial x_{k_{p}}^{0}}\right\vert_{*}
\end{equation}
The STTs have their own unique order-dependent dynamics. \cite{ParkScheeresSTTs} presents the dynamics for the STTs again element-by-element, up to fourth order. We compute them this way, getting the necessary Jacobian terms of Eq.~\eqref{STT3} (necessary for STT numerical integration) via automatic differentiation.

This STT representation is fundamentally redundant, e.g. the 3rd-order STT components associated with $\delta x_{1}^{0}\delta x_{2}^{0}\delta x_{1}^{0}$, $\delta x_{1}^{0}\delta x_{1}^{0}\delta x_{2}^{0}$, $\delta x_{2}^{0}\delta x_{1}^{0}\delta x_{1}^{0}$ are all equal yet appear individually within the STT structure. We can compute a fundamental solution from its associated STT terms as below:
\begin{equation}
\psi_{i,\mathcal{O}(x^{p})} = \frac{\mathcal{R}}{p!}\Phi_{i, k_{1}\ldots k_{p}}
\end{equation}
where $i = $1:$N$ denotes the $i$\textsuperscript{th} component (row) of the fundamental solution, $\mathcal{O}(x^{p})$ denotes some particular column of $\Psi_{j}$ associated with the order-$p$ STT component of interest, and $\mathcal{R}$ is the number of repeats of the relevant STT term -- for example, for $N = 2$, $j = 2$, and $\psi_{i,x_{1}x_{2}}$, we compute $\mathcal{R} = 2$ because we have two repeated STTs, $\Phi_{i,1\cdot2} = \Phi_{i,2\cdot1}$. Lastly recall that the ordering of the columns of $\Psi_{j}$ corresponds to the ordering of the corresponding monomials in $\bm{c}_{j}$.

\section{Computing the nonlinear fundamental solution expansions analytically}
Now we discuss a second method for constructing $\Psi_{j}(t)$. For the spacecraft relative motion problem there has been an enormous amount of past work to develop completely analytic nonlinear approximations of the solution of the satellite relative motion equations. These can be adopted directly to replace the numerically computed nonlinear functions used in this work. Here we provide a brief discussion about how such solutions can be derived. These ``solutions" approximate, via low-order series expansion in the initial states, the true behavior of the nonlinear differential equations of spacecraft relative motion given by \cite{Casotto}. Typically these are obtained by perturbation methods \citep{nayfeh2000perturbation,Hinch_Book}, most typically a straightforward expansion. Here we preview such a methodology briefly with a shortened discussion of how the second-order analytic STTs are derived \citep{Butcher2016_AAS}. In general, the state can be expanded in a perturbation series:
\begin{equation}
\label{pertSeries}
\bm{x}(\zeta) = \bm{x}_{1}(\zeta) + \varepsilon\bm{x}_{2}(\zeta) + \varepsilon^{2}\bm{x}_{3}(\zeta) + \ldots
\end{equation}
where $\varepsilon$ is a small parameter, $\zeta = nt$ is non-dimensional time, and typically the states are rendered in a non-dimensional form. The function $\bm{x}_{1}(\zeta)$ is the solution to the unperturbed linearized (CW) dynamics, and the subsequent functions correct for nonlinearities in the dynamics, and possibly also target orbit eccentricity.

First, the dimensionless CW equations are solved to obtain $\bm{x}_{1}(\zeta)$. Then, we perform a nonlinear expansion of the nonlinear equations of relative motion (see e.g. \cite{schaub} for Keplerian or \cite{Casotto} for non-Keplerian) up to a desired order in the states, then substitute the expansion of Eq.~\eqref{pertSeries} into the dynamics. The dynamics are then parsed order-by-order starting with the unperturbed linear problem $\mathcal{O}(\varepsilon^{0})$ and continuing to all other powers of $\varepsilon$. From the $\mathcal{O}(\varepsilon^{1})$ part, the model equations are rearranged into a form below separating all perturbations up to second-order in the states (constituting $F_{2}$) from the dimensionless Clohessy-Wiltshire (CW) ODE equations:

\begin{equation}
\label{FM_ch4_nondim}
\begin{split}
& \left( \begin{array}{c} x_{2}'' - 2y_{2}' - 3x_{2} \\ y_{2}'' + 2x_{2}' \\ z_{2}'' + z_{2} \end{array} \right) =  F_{2}\left(\zeta,x_{1},y_{1},z_{1},x_{1}',y_{1}',z_{1}' \right)
\end{split}
\end{equation}

The solutions to these equations can be found using symbolic software to evaluate the inverse Laplace transform:

\begin{equation}
\label{FM_ch4_pert}
\left( \begin{array}{c} x_{2}\left(\zeta\right) \\  y_{2}\left(\zeta\right) \\  z_{2}\left(\zeta\right) \end{array} \right) = \mathcal{L}^{-1}\left( \begin{bmatrix} s^{2} - 3 & -2s & 0 \\ 2s & s^{2} & 0 \\ 0 & 0 & s^{2} + 1 \end{bmatrix}^{-1} \mathcal{L}\left(F_{2}\left(\zeta\right)\right)\right)
\end{equation}

The inverted matrix is the transfer matrix of the CW system, and $F_{2}\left(\zeta\right)$ is obtained by substituting the normalized CW solution into all states in the right side of Eq. \eqref{FM_ch4_nondim}. The velocity states $x_{2}^{'}(\zeta)$, $y_{2}^{'}(\zeta)$, $z_{2}^{'}(\zeta)$ are then obtained by simple symbolic differentiation of the solutions above. 

For the analytic version of the third-order nonlinear solutions that we numerically computed to build $\Psi_{j}$, in local Cartesian coordinates, consult \cite{Butcher2016_AAS}. Note their solution is trigonometrically sorted for compactness instead of in terms of the monomials. \cite{ButcherBurnett2017} explores a similar approach in both Cartesian coordinates and curvilinear coordinates, which offer greater accuracy because they better accommodate the natural curvature of the orbit geometry. These solutions also apply corrections for nonzero target orbit eccentricity, which improves their fidelity in a true flight setting. Note that the eccentricity-linear part of the curvilinear solutions are derived from an equation which contains an error, which is addressed in the Appendix of \cite{willis2019secondorder}. That work also derives its own curvilinear solutions and provides a nice discussion of other analytic solutions in literature for this problem. Note that for long-duration problems, it becomes necessary to account for non-Keplerian perturbations such as $J_{2}$ \citep{Burnett2018_J2}. The appealing possibilities in use of these solutions in guidance applications  (particularly if they are analytic and time-explicit), as demonstrated by this paper, might breathe new life into the academic efforts for the derivation of nonlinear analytic solutions.

\section{Computations using Differential Algebra}
It is well-known that the STTs and Differential Algebra (DA) can be used to provide the same information about the nonlinear expansions in the vicinity of a reference. Thus DA constitutes a third method for building $\Psi_{j}(t)$. Indeed, we have tested this, computing $\Psi_{j}$ using Pyaudi \citep{audi_izzo_zenodo}, which offers speed advantages over the STTs. Similarly to the case of using STTs to build $\Psi_{j}$, in general $\Psi_{j}(t)$ can be constructed for discrete times $t_{i} \in [t_{0}, t_{1}, \ldots, t_{f}]$ by computing the Taylor map (TM) from $t_{0}$ to $t_{1}$, then $t_{0}$ to $t_{2}$, etc. for all times $t_{i} > t_{0}$, extracting the coefficients for every time, and ordering them appropriately to build a corresponding $\Psi_{j}(t_{i})$ $\forall$ $t_{i}$.  Here we contextualize some of our work within the DA framework as outlined in \cite{Berz_ParticleBeamMaps}. For this discussion we borrow and synthesize Berz' notation where needed. In lieu of a compressed summary, we defer the unfamiliar reader to their thorough introduction to DA.

First we note that our $K_{j}$-dimensional space $\mathcal{C}^{(N,j)}$ (facilitated by the $j$\textsuperscript{th}-order nonlinear expansion) is directly related to the vector space ${}_{j}D_{N}$. In particular, let $d_{k} = [x_{k}]$ denote the higher-order differential structure induced by a $j$\textsuperscript{th}-order expansion of the $k$\textsuperscript{th} element of coordinates $\bm{x} \in \mathbb{R}^{N}$. Then the Taylor expansion $T_{f}$ of some function $f$ can be expanded into the vector space ${}_{j}D_{N}$:
\begin{equation}
\label{DAc1}
[f] = [T_{f}] = \sum_{q_{1} + \ldots + q_{N} \leq j} \alpha_{q_{1},\ldots,q_{N}}\cdot d_{1}^{q_{1}}\ldots d_{N}^{q_{N}}
\end{equation}
Berz computes, via consideration of the tuples $q_{1},\ldots ,q_{N}$, that the dimensionality of ${}_{j}D_{N}$ is as follows:
\begin{equation}
\label{DAc2}
\text{dim}\left( {}_{j}D_{N} \right) = \begin{pmatrix} N + j \\ N \end{pmatrix}
\end{equation}
We observe the following similarity, where $\mathcal{C}^{(N,j)} \subseteq \mathbb{R}^{K_{j}}$, and $K_{j}$ is given by our Eq.~\eqref{WhatIsKj1}:
\begin{equation}
\label{DAc3}
\text{dim}\left( {}_{j}D_{N} \right)  = K_{j} + 1
\end{equation}
The difference here is that ${}_{j}D_{N}$ contains the ``Real" component (the reference point about which the nonlinear expansion occurs) by definition, but $\mathcal{C}^{(N,j)}$ does not, because our parameterization works only with deviations from a reference. 

\bibliographystyle{plainnat}
\end{document}